# CHARACTERISTIC CLASSES OF PROALGEBRAIC VARIETIES AND MOTIVIC MEASURES

SHOJI YOKURA

ABSTRACT. Michael Gromov has recently initiated what he calls "symbolic algebraic geometry", in which objects are proalgebraic varieties: a proalgebraic variety is by definition the projective limit of a projective system of algebraic varieties. In this paper we introduce characteristic classes of proalgebraic varieties, using Grothendieck transformations of Fulton–MacPherson's Bivariant Theory, modeled on the construction of MacPherson's Chern class transformation of proalgebraic varieties. We show that a proalgebraic version of the Euler–Poincaré characteristic with values in the Grothendieck ring is a generalization of the so-called motivic measure.

*Dedicated to Clint McCrory on the occasion of his 60th birthday*

## CONTENTS

## 1. INTRODUCTION

This work is motivated by Gromov's papers [Grom1, Grom2] and also by [Y7].

In [Grom1] Michael Gromov initiated what he calls "symbolic algebraic geometry" and in its Abstract he says " ... The paper intends to bring out relations between model theory, algebraic geometry, and symbolic dynamics." We hope that this present work would be an even tiny contribution to "symbolic algebraic geometry".

A *pro-algebraic variety* is defined to be a projective system of complex algebraic varieties and a *proalgebraic variety* is defined to be the projective limit of a pro-algebraic variety. Proalgebraic varieties are the main objects in [Grom1]. A pro-category was introduced by A. Grothendieck [Grot] and it was used to develop the Etale Homotopy Theory [AM] and Shape Theory (e.g., see [Bor], [Ed], [MS], etc.) and so on.

In [Grom1] M. Gromov investigated the *surjunctivity* [Got] , *i.e., being either surjective or non-injective,* in the category of proalgebraic varieties. The original or classical surjunctivity theorem is the so-called *Ax' Theorem*, saying that every regular selfmapping of a complex algebraic variety is surjunctive; thus if it is injective then it has to be surjective (cf, [Ax], [BBR], [Bo], [Kurd] , [New], [Par], etc.).

Our interest at the moment is not a further investigation concerning Ax-type theorems, but characteristic classes, in particular, Chern classes of proalgebraic varieties. A very

Partially supported by Grant-in-Aid for Scientific Research (No. 15540086, No.17540088), the Ministry of Education, Culture, Sports, Science and Technology (MEXT), Japan.

simple example of a proalgebraic variety is the Cartesian product $X^{\mathbb{N}}$ of an infinite countable copies of a complex algebraic variety $X$, which is one of the main objects treated in [Grom 1]. Then, what would be the Chern–Schwartz–MacPherson class of $X^{\mathbb{N}}$ ? In particular, what would be the *"Euler–Poincaré characteristic"* of $X^{\mathbb{N}}$ ? Our answers are that they are respectively the Chern–Schwartz–MacPherson class $c_*(X)$ and the Euler–Poincaré characteristic $\chi(X)$ in some sense, which will be clarified later. It is this very simple observation (which looked meaningless at the beginning) that led us to the present work, which naturally led us to motivic measures. The motivic measures/integrations have been actively studied by many people (e.g., see [Cr], [DL1], [DL2], [Kon], [Loo], [Ve] etc.).

In a general set-up we consider *bifunctors*. The bifunctors with which we deal are binfunctors $\mathcal{F} : \mathcal{C} \to \mathcal{A}$ from a category $\mathcal{C}$ to the category $\mathcal{A}$ of abelian groups, i.e., $\mathcal{F}$ is a pair $(\mathcal{F}_*, \mathcal{F}^*)$ of a covariant functor $\mathcal{F}_*$ and a contravariant functor $\mathcal{F}^*$ such that $\mathcal{F}_*(X) = \mathcal{F}^*(X)$ for any object $X$. Unless some confusion occurs, we just denote $\mathcal{F}(X)$ for $\mathcal{F}_*(X) = \mathcal{F}^*(X)$. A typical example is the constructible function functor $F(X)$. Furthermore we assume that for a final object $pt \in Obj(\mathcal{C})$, $\mathcal{F}(pt)$ is a commutative ring $\mathcal{R}$ with a unit. The morphism from an object $X$ to a final object $pt$ shall be denoted by $\pi_X : X \to pt$. Then the covariance of the bifunctor $\mathcal{F}$ induces the homomorphism $\mathcal{F}(\pi_X) : \mathcal{F}(X) \to \mathcal{F}(pt) = \mathcal{R}$, which shall be denoted by

$$\chi_{\mathcal{F}} : \mathcal{F}(X) \to \mathcal{R}$$

and called the *$\mathcal{F}$-characteristic*, just mimicking the Euler–Poincaré characteristic

$$\chi : F(X) \to \mathbb{Z}$$

in the case when $\mathcal{F} = F$.

Let $X_\infty = \varprojlim_{\lambda \in \Lambda} \left\{ X_\lambda, \pi_{\lambda\mu} : X_\mu \to X_\lambda \right\}$ be a proalgebraic variety and let $P = \left\{ p_{\lambda\mu} \right\}$ be a projective system of elements of $\mathcal{R}$ by the directed set $\Lambda$, i.e., a set such that $p_{\lambda\lambda} = 1$ (the unit) and $p_{\lambda\mu} \cdot p_{\mu\nu} = p_{\lambda\nu}$ $(\lambda < \mu < \nu)$. For each $\lambda \in \Lambda$ the subobject $\mathcal{F}_P^{\text{st}}(X_\lambda)$ of $\mathcal{F}(X_\lambda)$ is defined to be

$$\mathcal{F}_P^{\text{st}}(X_\lambda) := \left\{ \alpha_\lambda \in \mathcal{F}(X_\lambda) | \ \chi_{\mathcal{F}}\bigl(\pi_{\lambda\mu}{}^* \alpha_\lambda\bigr) = p_{\lambda\mu} \cdot \chi_{\mathcal{F}}(\alpha_\lambda) \text{ for any } \mu \text{ such that } \lambda < \mu \right\}.$$

The inductive limit $\varinjlim_\Lambda \left\{ \mathcal{F}_P^{\text{st}}(X_\lambda), \ \pi_{\lambda\mu}{}^* : \mathcal{F}_P^{\text{st}}(X_\lambda) \to \mathcal{F}_P^{\text{st}}(X_\mu) \ (\lambda < \mu) \right\}$ considered for a proalgebraic variety $X_\infty = \varprojlim_{\lambda \in \Lambda} X_\lambda$ is denoted by

$$\mathcal{F}_P^{\text{st.ind}}(X_\infty).$$

Our key observation, which is an application of standard facts on inductive systems and inductive limits, is the following:

**Theorem 1.1.** *(i) For a proalgebraic variety $X_\infty = \varprojlim_{\lambda \in \Lambda} \left\{ X_\lambda, \pi_{\lambda\mu} : X_\mu \to X_\lambda \right\}$ and a projective system $P = \left\{ p_{\lambda\mu} \right\}$ of non-zero elements of $\mathcal{R}$, we have the homomorphism*

$$\chi_{\mathcal{F}}^{\text{ind}} : \mathcal{F}_P^{\text{st.ind}}(X_\infty) \to \varinjlim_{\lambda \in \Lambda} \left\{ \times p_{\lambda\mu} : \mathcal{R} \to \mathcal{R} \right\},$$

*which is called the proalgebraic $\mathcal{F}$-characteristic homomorphism.*

*(ii) In the case when $\Lambda = \mathbb{N}$, for a proalgebraic variety $X_\infty = \varprojlim_{n \in \mathbb{N}} \left\{ X_n, \pi_{nm} : X_m \to X_n \right\}$ and a projective system $P = \{p_{nm}\}$ of non-zero elements of $\mathcal{R}$, the proalgebraic $\mathcal{F}$-characteristic homomorphism $\chi_{\mathcal{F}}^{\text{ind}} : \mathcal{F}_P^{\text{st.ind}}(X_\infty) \to \varinjlim_n \left\{ \times p_{nm} : \mathcal{R} \to \mathcal{R} \right\}$ is realized as the homomorphism*

$$\widetilde{\chi_{\mathcal{F}}^{\text{ind}}} : \mathcal{F}_P^{\text{st.ind}}(X_\infty) \to \mathcal{R}_P$$

*defined by*
$$\widetilde{\chi_{\mathcal{F}}^{\text{ind}}}\left([\alpha_n]\right) := \frac{\chi_{\mathcal{F}}(\alpha_n)}{p_{01} \cdot p_{12} \cdot p_{23} \cdots p_{(n-1)n}}.$$

Here $p_{01} := 1$ and $\mathcal{R}_P$ is the ring $\mathcal{R}_S$ of fractions of $\mathcal{R}$ with respect to the multiplicatively closed set $S$ consisting of all the finite products of powers of elements in $P$.

*(iii) In particular, in the case when the above projective system $P = \{p^s\}$ consists of powers of a non-zero element $p$, we get the homomorphism*
$$\widetilde{\chi_{\mathcal{F}}^{\text{ind}}} : \mathcal{F}_P^{\text{st.ind}}(X_\infty) \to \mathcal{R}\left[\frac{1}{p}\right]$$

*defined by*
$$\widetilde{\chi_{\mathcal{F}}^{\text{ind}}}\left([\alpha_n]\right) := \frac{\chi_{\mathcal{F}}(\alpha_n)}{p^{n-1}}.$$

Here $\mathcal{R}\left[\frac{1}{p}\right]$ is the ring of polynomials generated by $\{\frac{1}{p^s}\}$.

A typical example for Theorem 1.1 is the following. Let $X_\infty = \varprojlim_{n \in \mathbb{N}} \left\{X_n, \pi_{nm} : X_m \to X_n\right\}$ be a proalgebraic variety such that for each $n$ the structure morphism $\pi_{n(n+1)} : X_{n+1} \to X_n$ satisfies the condition that the Euler–Poincaré characteristics of the fibers of $\pi_{n(n+1)}$ are non-zero and constant; for example, $\pi_{n(n+1)} : X_{n+1} \to X_n$ is a locally trivial fiber bundle such that the Euler–Poincaré characteristic of the fiber is non-zero. Let us denote the constant Euler–Poincaré characteristic of the fibers of the morphism $\pi_{n(n+1)} : X_{n+1} \to X_n$ by $e_n$ and we set $e_0 := 1$. Then we get the canonical proalgebraic Euler–Poincaré characteristic homomorphism
$$\chi^{\text{ind}} : F^{\text{ind}}(X_\infty) \to \mathbb{Q}$$

described by
$$\chi^{\text{ind}}([\alpha_n]) = \frac{\chi(\alpha_n)}{e_0 \cdot e_1 \cdot e_2 \cdots e_{n-1}}.$$

In particular, if the Euler–Poincaré characteristics $e_n$ are all the same, say $e_n = e$ for each $n$, then the canonical proalgebraic Euler–Poincaré characteristic homomorphism $\chi^{\text{ind}} : F^{\text{ind}}(X_\infty) \to \mathbb{Q}$ is described by
$$\chi^{\text{ind}}([\alpha_n]) = \frac{\chi(\alpha_n)}{e^{n-1}}.$$

In this special case, the target ring $\mathbb{Q}$ can be replaced by the ring $\mathbb{Z}\left[\frac{1}{e}\right]$.

Let $K_0(\mathcal{V})$ be the Grothendieck ring of complex algebraic varieties. Then a "motivic" version of the Euler–Poincaré characteristic $\chi : F(X) \to \mathbb{Z}$ is the homomorphism $\Gamma : F(X) \to K_0(\mathcal{V})$ "tautologically" defined by $\Gamma(\sum_W a_W \mathbb{1}_W) := \sum_W a_W [W]$, where $\mathbb{1}_W$ denotes the characteristic function supported on a subvariety $W$ of $X$ and $a_W \in \mathbb{Z}$ and $[W] \in K_0(\mathcal{V})$. Then we get the following theorem, which is a generalization of the motivic measure:

**Theorem 1.2.** *(i) For a proalgebraic variety $X_\infty = \varprojlim_{\lambda \in \Lambda} \left\{X_\lambda, \pi_{\lambda\mu} : X_\mu \to X_\lambda\right\}$ and a projective system $G = \{\gamma_{\lambda\mu}\}$ of non-zero Grothendieck classes, we get the proalgebraic Grothendieck class homomorphism*
$$\Gamma^{\text{ind}} : F_G^{\text{st.ind}}(X_\infty) \to \varinjlim_{\lambda \in \Lambda} \left\{\times \gamma_{\lambda\mu} : K_0(\mathcal{V}) \to K_0(\mathcal{V})\right\}.$$

*(ii) In the case when $\Lambda = \mathbb{N}$, for a proalgebraic variety $X_\infty = \varprojlim_{n \in \mathbb{N}} \left\{X_n, \pi_{nm} : X_m \to X_n\right\}$ and a projective system $G = \{\gamma_{n,m}\}$ of non-zero Grothendieck classes, we* *have the following canonical proalgebraic Grothendieck class homomorphism*

$$\widetilde{\Gamma^{\text{ind}}} : F_G^{\text{st.ind}}(X_\infty) \to K_0(\mathcal{V})_G$$

*which is defined by*

$$\widetilde{\Gamma^{\text{ind}}}\left([\alpha_n]\right) := \frac{\Gamma(\alpha_n)}{\gamma_{01} \cdot \gamma_{12} \cdot \gamma_{23} \cdots \gamma_{(n-1)n}}.$$

*Here we set $\gamma_{01} := \mathbb{1}$ and $K_0(\mathcal{V})_G$ is the ring of fractions of $K_0(\mathcal{V})$ with respect to the multiplicatively closed set consisting of finite products of powers of elements of $G$.*

(iii) *Let $X_\infty = \varprojlim_{n \in \mathbb{N}} \left\{ X_n, \pi_{nm} : X_m \to X_n \right\}$ be a proalgebraic variety such that each structure morphism $\pi_{n(n+1)} : X_{n+1} \to X_n$ satisfies the condition that for each $n$ there exists a $\gamma_n \in K_0(\mathcal{V})$ such that $\pi_{n(n+1)}^{-1}(S_n) = \gamma_n \cdot [S_n]$ for any constructible set $S_n \subset X_n$; for example, $\pi_{n(n+1)} : X_{n+1} \to X_n$ is a Zariski locally trivial fiber bundle with fiber variety being $F_n$ (in which case $\gamma_n = [F_n] \in K_0(\mathcal{V})$). Then the canonical proalgebraic Grothendieck class homomorphism*

$$\Gamma^{\text{ind}} : F^{\text{ind}}(X_\infty) \to K_0(\mathcal{V})_G$$

*is described by*

$$\Gamma^{\text{ind}}([\alpha_n]) = \frac{\Gamma(\alpha_n)}{\gamma_0 \cdot \gamma_1 \cdot \gamma_2 \cdots \gamma_{n-1}}.$$

*Here $\gamma_0 := \mathbb{1}$ and $K_0(\mathcal{V})_G$ is the ring of fractions of $K_0(\mathcal{V})$ with respect to the multiplicatively closed set consisting of finite products of powers of $\gamma_m$ ($m = 1, 2, 3 \cdots$).*

(iv) *In particular, if $\gamma_n$ are all the same, say $\gamma_n = \gamma$ for any $n$, then the canonical proalgebraic Grothendieck class homomorphism*

$$\Gamma^{\text{ind}} : F^{\text{ind}}(X_\infty) \to K_0(\mathcal{V})_G$$

*is described by*

$$\Gamma^{\text{ind}}([\alpha_n]) = \frac{\Gamma(\alpha_n)}{\gamma^{n-1}}.$$

*In this special case the quotient ring $K_0(\mathcal{V})_G$ shall be simply denoted by $K_0(\mathcal{V})_\gamma$.*

In passing, here it may be instructive to remark the following from symbolic dynamics. The standard metric to be considered on the sequence space

$$2^{\mathbb{N}} := \{(x_0, x_1, x_2, \cdots, x_n, \cdots) \mid x_i \in \{0, 1\}\}$$

is: for $a = (a_0, a_1, a_2, \cdots), b = (b_0, b_1, b_2, \cdots)$,

$$d(a, b) := \sum_{n=0}^{\infty} \frac{|a_n - b_n|}{2^n}.$$

This metric now looks very natural, since the sequence space $2^{\mathbb{N}}$ is proalgebraic and the Euler–Poincaré characteristic $\chi(\{0, 1\}) = 2$. More generally, for the space $k^{\mathbb{N}} := \{(x_0, x_1, x_2, \cdots, x_n, \cdots) \mid k_i \in \{0, 1, 2, \cdots, k-1\}\}$, the standard metric on this sequence space is the same as above. However, from our viewpoint, we can consider also the following metric:

$$d_k(a, b) := \sum_{n=0}^{\infty} \frac{|a_n - b_n|}{k^n},$$

because $\chi(\{0, 1, \cdots, k-1\}) = k$

When we extend MacPherson's Chern class transformation [Mac1] to a category of proalgebraic varieties, we appeal to the Bivariant Theory, which was introduced by W. Fulton and R. MacPherson [FM] as a theory unifying both covariant and contravariant functors in order to apply to singular spaces. A generalized formulation of characteristic classes of proalgebraic varieties is the following

**Theorem 1.3.** *(i) Let $\gamma : \mathbb{B} \to \mathbb{B}'$ be a Grothendieck transformation between two bivariant theories $\mathbb{B}, \mathbb{B}' : \mathcal{C} \to \mathcal{A}$ and let $\{(\pi_{\lambda\mu}; b_{\lambda\mu}) : X_\mu \to X_\lambda\}$ be a projective system of bivariant-class-equipped morphisms. Then we get the following pro-version of the natural transformation $\gamma_* : \mathbb{B}_* \to \mathbb{B}'_*$:*

$$\gamma_*^{\text{ind}} : \mathbb{B}_*^{\text{ind}}\left(X_\infty; \{b_{\lambda\mu}\}\right) \to \mathbb{B}'^{\text{ind}}_*\left(X_\infty; \{\gamma(b_{\lambda\mu})\}\right).$$

*(ii) Let $\{f_\lambda : Y_\lambda \to X_\lambda\}$ be a fiber-square pro-morphism between two projective systems $\{(\rho_{\lambda\mu}; d_{\lambda\mu}) : Y_\mu \to Y_\lambda\}$ and $\{(\pi_{\lambda\mu}; b_{\lambda\mu}) : X_\mu \to X_\lambda\}$ of bivariant-class-equipped morphisms such that $d_{\lambda\mu} = f_\lambda^\star b_{\lambda\mu}$. Then we have the following commutative diagram:*

$$\begin{array}{ccc}
\mathbb{B}_*^{\text{ind}}(Y_\infty; \{d_{\lambda\mu}\}) & \xrightarrow{\gamma_*^{\text{ind}}} & \mathbb{B}'^{\text{ind}}_*(Y_\infty; \{\gamma(d_{\lambda\mu})\}) \\
f_{\infty*} \downarrow & & \downarrow f_{\infty*} \\
\mathbb{B}_*^{\text{ind}}(X_\infty; \{b_{\lambda\mu}\}) & \xrightarrow[\gamma_*^{\text{ind}}]{} & \mathbb{B}'^{\text{ind}}_*(X_\infty; \{\gamma(b_{\lambda\mu})\}).
\end{array}$$

*(iii) Let $\mathbb{B}_*(pt) = \mathbb{B}'_*(pt)$ be a commutative ring $\mathcal{R}$ with a unit and we assume that the homomorphism $\gamma : \mathbb{B}_*(pt) \to \mathbb{B}'_*(pt)$ is the identity. Let $P = \{p_{\lambda\mu}\}$ be a projective system of non-zero elements $p_{\lambda\mu} \in \mathcal{R}$. Then we get the commutative diagram*

$$\mathbb{B}_{*,P}^{\text{st.ind}}\left(X_\infty; \{b_{\lambda\mu}\}\right) \xrightarrow{\gamma_*^{\text{ind}}} \mathbb{B}'^{\text{st.ind}}_{*,P}\left(X_\infty; \{\gamma(b_{\lambda\mu})\}\right)$$

with $\chi_{\mathbb{B}_*}^{\text{ind}}$ and $\chi_{\mathbb{B}'_*}^{\text{ind}}$ mapping to $\varinjlim_{\lambda \in \Lambda}\left\{\times p_{\lambda\mu} : \mathcal{R} \to \mathcal{R}\right\}.$

If we apply this theorem to Brasselet's bivariant Chern class [Br] or [BSY1], we get a proalgebraic version of MacPherson's Chern class transformation $c_* : F \to H_*$.

The relative version of the above Grothendieck ring $K_0(\mathcal{V})$ is the relative Grothendieck ring $K_0(\mathcal{V}/X)$ of complex algebraic varieties over a variety $X$, which is a bifunctor, and there is a canonical homomorphism $e : K_0(\mathcal{V}/X) \to F(X)$ defined by $e([Y \xrightarrow{f} X]) := f_* \mathbb{1}_Y$. This is a natural transformation and in §6 we will show that this natural transformation is unique in a sense.

*Acknowledgements.* I would like to thank Jörg Schürmann and Willem Veys for their valuable comments and suggestions on an earlier version (2003). A part of the paper was written during the author's visits at the Erwin Schrödinger International Institute of Mathematical Physics (ESI), Vienna, in August 2004 and 2005. I would like to thank the staff of the ESI and Professor Peter Michor for providing a nice atmosphere in which to work. Finally I would like to give my sincere gratitude to Clint McCrory for having introduced me to characteristic classes of singular spaces in my graduate student days.

## 2. Constructible functions and MacPherson's Chern class transformation

A constructible set in an algebraic variety is one obtained from the subvarieties by taking finitely many of the operations of intersection $\cap$, union $\cup$, subtraction $-$. The collection of such sets is sometimes called the Boolean algebra of $X$ generated by the subvarieties of $X$. A constructible function on a variety is an integer-valued function for which the variety has a finite stratification into constructible sets such that the function is constant on each constructible set. The abelian group of all constructible functions on a variety $X$ is denoted by $F(X)$. Equivalently we can describe the group $F(X)$ as follows. For a subvariety $W$ of a given variety $X$, $\mathbb{1}_W$ denotes the characteristic function supported on the subvariety $W$, i.e., $\mathbb{1}_W(x) = 1$ for $x \in W$ and $\mathbb{1}_W(x) = 0$ for $x \notin W$. Then $F(X)$ consists of all finite linear combinations of such characteristic functions supported on subvarieties with integer coefficients.

If we define the pullback $f^* : F(Y) \to F(X)$ by the usual functional pullback, i.e., $f^*\alpha := \alpha \circ f$ for a constructilbe function $\alpha \in F(Y)$, then the assignment $X \longmapsto F(X)$ is *a contravariant functor*; for morphisms $f : X \to Y$ and $g : Y \to Z$ we have $(g \circ f)^* = f^* \circ g^*$.

It turns out that the assignment $X \longmapsto F(X)$ is also *a covariant functor*:

**Proposition 2.1.** *([Mac1, Proposition 1]) For a morphism $f : X \to Y$, we define the pushforward*

$$f_* : F(X) \to F(Y) \quad by \quad f_*(\mathbb{1}_W)(p) = \chi(f^{-1}(p) \cap W).$$

*Then it is well-defined and it is covariantly functorial, i.e., for morphisms $f : X \to Y$ and $g : Y \to Z$ we have $(g \circ f)_* = g_* \circ f_*$.*

Let us define the Euler–Poincaré characteristic homomorphism $\chi : F(X) \to \mathbb{Z}$ by

$$\chi(\alpha) := \sum_{n \in \mathbb{Z}} n\chi(\alpha^{-1}(n)).$$

Then for a morphism $p : X \to pt$ to a point $pt$, the pushforward $p_* : F(X) \to F(pt) = \mathbb{Z}$ is nothing but the above $\chi : F(X) \to \mathbb{Z}$. So, if we consider the morphism $g : Y \to pt$ to a point $pt$ in the above equality $(g \circ f)_* = g_* \circ f_*$, we get the commutative diagram:

$$\begin{array}{ccc} F(X) & & \\ {\scriptstyle f_*}\downarrow & \searrow^{\chi} & \\ & & \mathbb{Z}. \\ & \nearrow_{\chi} & \\ F(Y) & & \end{array}$$

In fact, the commutativity of this diagram follows from the definition of the pushforward $f_* : F(X) \to F(Y)$ and the stratification theory; it is implicit in the proof of the above proposition.

For more details on constructible functions and, in particular, for comparison with standard Grothendieck operations on constructible sheaves, see [Dim], [KS], [Scha] and [Schü3].

What P. Deligne and A. Grothendieck conjectured and later R. MacPherson affirmatively solved is the following:

**Theorem 2.2.** *([Mac1, Theorem 1]) There exists a unique natural transformation from the covariant constructible function functor to the Borel–Moore homology covariant functor*

$$c_* : F \to H_*$$

*such that for a nonsingular variety $X$ the value of the characteristic function $\mathbb{1}_X$ is the Poincaré dual of the total Chern cohomology class:*

$$c_*(\mathbb{1}_X) = c(TX) \cap [X]$$

*where $TX$ is the tangent bundle of $X$.*

The formulation of the natural transformation $c_* : F \to H_*$ was motivated by that of the Stiefel–Whitney classes in the real case due to D. Sullivan [Sull] (also see [Fu-Mc]).

The above theorem is an answer for the question of whether or not there exists (uniquely) a homomorphism $? : F(X) \to H_*(X)$ such that the following diagram commutes

$$\begin{array}{ccc} F(X) & \xrightarrow{?} & H_*(X) \\ & \chi \searrow \swarrow \int_X & \\ & \mathbb{Z} & \end{array}$$

and such that it is functorial, i.e., the following diagram commutes:

$$\begin{array}{ccc} F(X) & \xrightarrow{?} & H_*(X) \\ f_* \downarrow & & \downarrow f_* \\ F(Y) & \xrightarrow{?} & H_*(Y). \end{array}$$

Here $\int_X : H_*(X) \to \mathbb{Z}$ is the integration or equal to $(\pi_X)_* : H_*(X) \to H_*(pt) = \mathbb{Z}$ with $\pi_X : X \to pt$ being the map to a point $pt$. It is obviously a Riemann–Roch type question for Chern classes just like Grothendieck extended the celebrated Hirzebruch–Riemann–Roch theorem to the so-called Grothendieck–Riemann–Roch theorem.

The key ingredients of MacPherson's proof of the above theorem are *local Euler obstruction*, *Chern–Mather class* and the *graph construction*.

The local Euler obstruction $\operatorname{Eu}_W$ supported on a subvariety $W$ is a certain constructible function supported on $W$, defined in an obstruction-theoretical way, which is identical to the charactristic function $\mathbb{1}_W$ *off* the singularities of $W$. Thus the values of $\operatorname{Eu}_W$ reflect singularities. The Chern–Mather class of $W$, denoted by $c^M(W)$, is defined by the Nash blow-up $\nu : \widehat{W} \to W$. To be a little more precise, if we let $\widehat{TW}$ be the tautological Nash tangent bundle on the Nash blow-up $\widehat{W}$, then it is defined by

$$c^M(W) = \nu_*(c(\widehat{TW}) \cap [\widehat{W}]).$$

One can see, by the induction on dimensions, that the abelian group $F(X)$ of constructible functions on $X$ is generated also by the local Euler obstructions supported on subvarieties. And it is an ingenious insight of MacPherson that the assignment $\operatorname{Eu}_W \longmapsto c^M(W)$, *instead of* the assignment $\mathbb{1}_W \longmapsto c^M(W)$, from $F(X)$ to the homology group $H_*(X)$ gives rise to *the looked for natural transformation* $c_* : F \to H_*$, which is proved by the graph construction.

For a characteristic function $\mathbb{1}_X$ the value $c_*(\mathbb{1}_X)$ is denoted simply by $c_*(X)$ and called the *MacPherson–Chern* (or *Chern–MacPherson*) class of the variety $X$. In particular, by considering the morphism to a point, the degree of the zero-dimensional component of $c_*(X)$ is nothing but the Euler–Poincaré charcateristic $\chi(X)$ of the variety $X$. In this sense the Chern–Schwartz–MacPherson class of a variety is a "higher homological class version" of the Euler–Poincaré characteristic of the given variety, just like the classical Chern class of a complex manifold is a "higher cohomological class version" of the the Euler–Poincaré characteristic of the manifold. The degree of the zero-dimensional component of the Chern–Mather class $c^M(X)$ is provisionally called the *Euler–Poincaré–Mather characteristic* of the variety $X$ and denoted by $\chi^M(X)$. Thus we have that $\chi(X) = \chi^M(X) +$ some numbers reflecting singularities.

Later J.-P. Brasselet and M.-H. Schwartz [BS] showed that for a variety $X$ embedded in a manifold $M$ the MacPherson–Chern class $c_*(X)$ of $X$ is isomorphic to the Schwartz class $c^{\operatorname{Sch}}(X) \in H^*(M, M \setminus X)$ of $X$ by the Alexander duality isomorphism $H^*(M, M \setminus X) \cong H_*(X)$. Thus the above natural transformation $c_* : F \to H_*$ is called the *MacPherson's Chern class transformation*.

**Remark 2.3.** In [STV] J. Seade, M. Tibar and A. Verjovsky introduced the notion of *global Euler obstruction* and showed that it is equal to the Euler–Poincaré–Mather characteristic.

## 3. Indconstructible functions, Euler–Poincaré characteristics of proalgebraic varieties and their generalizations

Let $I$ be a directed set and let $\mathcal{C}$ be a given category. Then a projective system is, by definition, a system $\{X_i, \pi_{ii'} : X_{i'} \to X_i (i < i'), I\}$ consisting of objects $X_i \in \text{Obj}(\mathcal{C})$, morphisms $\pi_{ii'} : X_{i'} \to X_i \in \text{Mor}(\mathcal{C})$ for each $i < i'$ and the index set $I$. The object $X_i$ is called a *term* and the morphism $\pi_{ii'} : X_{i'} \to X_i$ a *bonding morphism* or *structure morphism* ([MS]). The projective system $\{X_i, \pi_{ii'} : X_{i'} \to X_i (i < i'), I\}$ is sometimes simply denoted by $\{X_i\}_{i \in I}$.

Given a category $\mathcal{C}$, Pro-$\mathcal{C}$ is the category whose objects are projective systems $X = \{X_i\}_{i \in I}$ in $\mathcal{C}$ and whose set of morphisms from $X = \{X_i\}_{i \in I}$ to $Y = \{Y_j\}_{j \in J}$ is

$$\text{Pro-}\mathcal{C}(X, Y) := \varprojlim_J (\varinjlim_I \mathcal{C}(X_i, Y_j)).$$

This definition is not crystal clear, but a more down-to-earth definition is the following (e.g., see [Fox] or [MS]): A morphism $f : X \to Y$ consists of a map $\theta : J \to I$ (not necessarily order preserving) and morphisms $f_j : X_{\theta(j)} \to Y_j$ for each $j \in J$, subject to the condition that if $j < j'$ in $J$ then for some $i \in I$ such that $i > \theta(j)$ and $i > \theta(j')$, the following diagram commutes

$$\begin{array}{ccc}
& X_i & \\
\pi_{\theta(j')i} \swarrow & & \searrow \pi_{\theta(j)i} \\
X_{\theta(j')} & & X_{\theta(j)} \\
f_{j'} \downarrow & & \downarrow f_j \\
Y_{j'} & \xrightarrow{\rho_{jj'}} & Y_j
\end{array}$$

Given a projective system $X = \{X_i\}_{i \in I} \in \text{Pro-}\mathcal{C}$, the projective limit $X_\infty := \varprojlim X_i$ may not belong to the source category $\mathcal{C}$. For a certain sufficient condition for the existence of the projective limit in the category $\mathcal{C}$, see [MS] for example.

An object in Pro-$\mathcal{C}$ is called a *pro-object*. A projective system of algebraic varieties is called a *pro-algebraic variety* and its projective limit is called a *proalgebraic variety*, which may not be an algebraic variety but simply a topological space.

A pro-morphism between two pro-objects is quite complicated, as remarked above. However, it follows from [MS] that the pro-morphism can be described more naturally as a so-called *level preserving pro-morphism*. Suppose that we have two pro-algebraic varieties $X = \{X_\gamma\}_{\gamma \in \Gamma}$ and $Y = \{Y_\lambda\}_{\lambda \in \Lambda}$. Then a pro-algebraic morphism $\Phi = \{f_\lambda\}_{\lambda \in \Lambda} : X \to Y$ is described as follows: there is an order-preserving map $\xi : \Lambda \to \Gamma$, i.e., $\xi(\lambda) < \xi(\mu)$ for $\lambda < \mu$, and for each $\lambda \in \Lambda$ there is a morphism $f_\lambda : X_{\xi(\lambda)} \to Y_\lambda$ such that for $\lambda < \mu$ the following diagram commutes:

$$\begin{array}{ccc}
X_{\xi(\mu)} & \xrightarrow{f_\mu} & Y_\mu \\
\pi_{\xi(\lambda)\xi(\mu)} \downarrow & & \downarrow \rho_{\lambda\mu} \\
X_{\xi(\lambda)} & \xrightarrow{f_\lambda} & Y_\lambda,
\end{array}$$

Then, the projective limit of the system $\{f_\lambda\}$ is a morphism from the proalgebraic variety $X_\infty = \varprojlim_{\lambda \in \Lambda} X_\lambda$ to the proalgebraic variety $Y_\infty = \varprojlim_{\gamma \in \Gamma} Y_\gamma$. It is called a proalgebraic morphism and denoted by $f_\infty : X_\infty \to Y_\infty$.

From now on, to make the presentation simpler, we assume that a pro-morphism (pro-morphism, resp.) is (the projective limit of, resp.) a projective system of morphisms of varieties with the same directed set and that the above order-preserving map $\xi : \Lambda \to \Lambda$ is the identity.

Let $T : \mathcal{C} \to \mathcal{D}$ be a covariant functor between two categories $\mathcal{C}, \mathcal{D}$. Obviously the covariant functor $T$ extends to a covariant pro-functor
$$\text{Pro-}T : \text{Pro-}\mathcal{C} \to \text{Pro-}\mathcal{D}$$
defined by $\text{Pro-}T(\{X_\lambda\}_{\lambda \in \Lambda}) := \{T(X_\lambda)\}_{\lambda \in \Lambda}$. Let $T_1, T_2 : \mathcal{C} \to \mathcal{D}$ be two covariant functors and $N : T_1 \to T_2$ be a natural transformation between the two functors $T_1$ and $T_2$. Then the natural transformation $N : T_1 \to T_2$ extends to a natural pro-transformation
$$\text{Pro-}N : \text{Pro-}T_1 \to \text{Pro-}T_2.$$

Thus a pro-algebraic version of MacPherson's Chern class transformation is straightforward, i.e., we have
$$\text{Pro-}c_* : \text{Pro-}F \to \text{Pro-}H_*.$$
In this case, the characteristic pro-function $\mathbb{1}_X$ of the pro-algebraic variety $X = \{X_\lambda\}_{\lambda \in \Lambda}$ should be simply $\mathbb{1}_X := \{\mathbb{1}_{X_\lambda}\}_{\lambda \in \Lambda}$ and thus the pro-version of MacPherson's Chern class transformation of the pro-algebraic variety $X = \{X_\lambda\}_{\lambda \in \Lambda}$ is simply $\text{Pro-}c_*(X) = \{c_*(X_\lambda)\}_{\lambda \in \Lambda}$.

Furthermore, taking the projective limit of the above projective system of natural transformations $\text{Pro-}c_* : \text{Pro-}F \to \text{Pro-}H_*$ gives rise to a natural transformation
$$\varprojlim_{\lambda \in \Lambda} c_* : \varprojlim_{\lambda \in \Lambda} F(X_\lambda) \to \varprojlim_{\lambda \in \Lambda} H_*(X_\lambda).$$

**Remark 3.1.** In Etale Homotopy Theory [AM] and Shape Theory (e.g., see [Bor], [Ed], [MS]) one stays in the pro-category and does not consider limits and colimits, because doing so throw away some geometric informations.

The covariance of $F$ gives rise to the projective limit $\varprojlim_{\lambda \in \Lambda} F(X_\lambda)$ on one hand. On the other hand, since $F$ is also a contravariant functor, it is reasonable to define the following:

**Definition 3.2.** For a proalgebraic variety $X_\infty = \varprojlim_{\lambda \in \Lambda} \{X_\lambda, \pi_{\lambda\mu} : X_\mu \to X_\lambda\}$, the inductive limit of the inductive system $\{F(X_\lambda), \pi_{\lambda\mu}{}^* : F(X_\lambda) \to F(X_\mu)(\lambda < \mu)\}$ is denoted by $F^{\text{ind}}(X_\infty)$, i.e.,
$$F^{\text{ind}}(X_\infty) := \varinjlim_{\lambda \in \Lambda} F(X_\lambda) = \bigcup_\mu \rho^\mu\bigl(F(X_\mu)\bigr)$$
where $\rho^\mu : F(X_\mu) \to \varinjlim_{\lambda \in \Lambda} F(X_\lambda)$ is the homomorphism sending $\alpha_\mu$ to its equivalence class $[\alpha_\mu]$ of $\alpha_\mu$. An element of the group $F^{\text{ind}}(X_\infty)$ is called a *indconstructible function* on the proalgebraic variety $X_\infty$.

**Remark 3.3.** (i) $F^{\text{ind}}(X_\infty)$ does depend on the given projective system $\mathcal{S} = \{X_\lambda, \pi_{\lambda\mu} : X_\mu \to X_\lambda \ (\lambda < \mu)\}$; so in this sense it should be denoted by something like $F^{\text{ind}}_{\mathcal{S}}(X_\infty)$ with the reference to the system $\mathcal{S}$, but for the sake of simplicity we drop the subscript $\mathcal{S}$.

(ii) In an earlier version (math.AG/0407237) the above inductive limit was denoted by $F^{\text{pro}}(X_\infty)$ and an element of it was called a *proconstructible function*. But in this revised version we use the qualifier *indconstructible* becaus it is defined via the *inductive* limits, and the term "proconstructible function" will be used in a different context in a later section.

Definition 3.2 can be used for any contravariant functor. Namely, if $\mathcal{F} : \mathcal{C} \to \mathcal{A}$ is a contravariant functor and $\{X_\lambda, \pi_{\lambda\mu} : X_\mu \to X_\lambda (\lambda < \mu)\}$ is a projective system in $\mathcal{C}$, then for the projective limit $X_\infty = \varprojlim_{\lambda \in \Lambda} X_\lambda$, which itself *may not belong to* the category $\mathcal{C}$, we can define

$$\mathcal{F}^{\mathrm{ind}}(X_\infty) := \varinjlim_{\lambda \in \Lambda} \left\{ \mathcal{F}(X_\lambda), {\pi_{\lambda\mu}}^* : \mathcal{F}(X_\lambda) \to \mathcal{F}(X_\mu)(\lambda < \mu) \right\},$$

which also *may not belong to* the category $\mathcal{A}$.

In order to go further, we require $\mathcal{F} : \mathcal{C} \to \mathcal{A}$ to be a *bifunctor* from a category $\mathcal{C}$ to the category $\mathcal{A}$ of abelian groups, i.e., $\mathcal{F}$ is a pair $(\mathcal{F}_*, \mathcal{F}^*)$ of a covariant functor $\mathcal{F}_*$ and a contravariant functor $\mathcal{F}^*$ such that $\mathcal{F}_*(X) = \mathcal{F}^*(X)$ for any object $X$. Unless some confusion occurs, we just denote $\mathcal{F}(X)$ for $\mathcal{F}_*(X) = \mathcal{F}^*(X)$. Furthermore we assume that for a final object $pt \in Obj(\mathcal{C})$, $\mathcal{F}(pt)$ is a commutative ring $\mathcal{R}$ with a unit. The morphism from an object $X$ to a final object $pt$ shall be denoted by $\pi_X : X \to pt$. Then the covariance of the bifunctor $\mathcal{F}$ induces the homomorphism $\mathcal{F}(\pi_X) : \mathcal{F}(X) \to \mathcal{F}(pt) = \mathcal{R}$, which shall be denoted by

$$\chi_\mathcal{F} : \mathcal{F}(X) \to \mathcal{R}$$

and called the $\mathcal{F}$-*characteristic*, just mimicking the Euler–Poincaré characteristic $\chi : F(X) \to \mathbb{Z}$ in the case when $\mathcal{F} = F$. Then furthermore the covariance of $\mathcal{F}$ implies that for a morphism $f : X \to Y$ in $Mor(\mathcal{C})$ we get the commutative diagram

$$\begin{array}{ccc} \mathcal{F}(X) & & \\ & \searrow^{\chi_\mathcal{F}} & \\ f_* \downarrow & & \mathcal{R}. \\ & \nearrow_{\chi_\mathcal{F}} & \\ \mathcal{F}(Y) & & \end{array}$$

In Theorem 3.4 below we do not need the commutativity of this diagram or the commtativity of $\chi_\mathcal{F}$ with the pushforward, but here we mention this as an analogy of the case of the constructible function functor $F$. We come back to this commutative diagram later in this section and in §5 when we discuss functorialities of $\mathcal{F}^{\mathrm{ind}}$ and so forth.

Let $X_\infty = \varprojlim_{\lambda \in \Lambda} \left\{ X_\lambda, \pi_{\lambda\mu} : X_\mu \to X_\lambda \right\}$ be a proalgebraic variety and let $P = \{p_{\lambda\mu}\}$ be a projective system of elements of $\mathcal{R}$ by the directed set $\Lambda$, i.e., a set such that

$$p_{\lambda\lambda} = 1 \quad \text{(the unit)} \quad \text{and} \quad p_{\lambda\mu} \cdot p_{\mu\nu} = p_{\lambda\nu} \quad (\lambda < \mu < \nu).$$

For each $\lambda \in \Lambda$ we define the following subobject of $\mathcal{F}(X_\lambda)$:

$$\mathcal{F}^{\mathrm{st}}_P(X_\lambda) := \left\{ \alpha_\lambda \in \mathcal{F}(X_\lambda) \mid \chi_\mathcal{F}\bigl({\pi_{\lambda\mu}}^* \alpha_\lambda\bigr) = p_{\lambda\mu} \cdot \chi_\mathcal{F}(\alpha_\lambda) \text{ for any } \mu \text{ such that } \lambda < \mu \right\}.$$

For each $\lambda \in \Lambda$, an element of $\mathcal{F}^{\mathrm{st}}_P(X_\lambda)$ is called *a $\chi_\mathcal{F}$-stable object of $\mathcal{A}$ with respect to the projective system $P$*. Then it is easy to see that for each structure morphism $\pi_{\lambda\mu} : X_\mu \to X_\lambda$ the pullback homomorphism ${\pi_{\lambda\mu}}^* : \mathcal{F}(X_\lambda) \to \mathcal{F}(X_\mu)$ preserves $\chi_\mathcal{F}$-stable objects with respect to the projective system $\{p_{\lambda\mu}\}$, namely it induces the homomorphism (using the same symbol):

$$ {\pi_{\lambda\mu}}^* : \mathcal{F}^{\mathrm{st}}_P(X_\lambda) \to \mathcal{F}^{\mathrm{st}}_P(X_\mu)$$

which implies that we get the inductive system

$$\left\{ \mathcal{F}^{\mathrm{st}}_P(X_\lambda), \quad {\pi_{\lambda\mu}}^* : \mathcal{F}^{\mathrm{st}}_P(X_\lambda) \to \mathcal{F}^{\mathrm{st}}_P(X_\mu) \quad (\lambda < \mu) \right\}.$$

Then for a proalgebraic variety $X_\infty = \varprojlim_{\lambda \in \Lambda} X_\lambda$ we consider the inductive limit of the above inductive system and it shall be denoted by

$$\mathcal{F}_P^{\text{st.ind}}(X_\infty)$$

and an element of this inductive limit shall be called a $\chi_\mathcal{F}$-*stable indobject* of $\mathcal{A}$ *on the proalgebraic variety* $X_\infty$ *with respect to the projective system* $P$. We see that this can be also directly defined as follows:

$$\left\{ [\alpha_\lambda] \in \mathcal{F}^{\text{ind}}(X_\infty) \quad | \quad \chi_\mathcal{F}(\pi_{\lambda\mu}{}^* \alpha_\lambda) = p_{\lambda\mu} \cdot \chi_\mathcal{F}(\alpha_\lambda) \quad (\lambda < \mu) \right\}.$$

The following is an application of standard facts on indutive systems and inductive limits, but nevertheless it is a key and important observation for the rest of the paper, in particular in connection to motivic measures, so it is stated as a theorem.

**Theorem 3.4.** *(i) For a proalgebraic variety* $X_\infty = \varprojlim_{\lambda \in \Lambda} \left\{ X_\lambda, \pi_{\lambda\mu} : X_\mu \to X_\lambda \right\}$ *and a projective system* $P = \{p_{\lambda\mu}\}$ *of non-zero elements of* $\mathcal{R}$, *we have the homomorphism*

$$\chi_\mathcal{F}^{\text{ind}} : \mathcal{F}_P^{\text{st.ind}}(X_\infty) \to \varinjlim_{\lambda \in \Lambda} \left\{ \times p_{\lambda\mu} : \mathcal{R} \to \mathcal{R} \right\},$$

*which is called the proalgebraic* $\mathcal{F}$-*characteristic homomorphism.*

*(ii) In the case when* $\Lambda = \mathbb{N}$, *for a proalgebraic variety* $X_\infty = \varprojlim_{n \in \mathbb{N}} \left\{ X_n, \pi_{nm} : X_m \to X_n \right\}$ *and a projective system* $P = \{p_{nm}\}$ *of non-zero elements of* $\mathcal{R}$, *a proalgebraic* $\mathcal{F}$-*characteristic homomorphism* $\chi_\mathcal{F}^{\text{ind}} : \mathcal{F}_P^{\text{st.ind}}(X_\infty) \to \varinjlim_n \left\{ \times p_{nm} : \mathcal{R} \to \mathcal{R} \right\}$ *is realized as the homomorphism*

$$\widetilde{\chi_\mathcal{F}^{\text{ind}}} : \mathcal{F}_P^{\text{st.ind}}(X_\infty) \to \mathcal{R}_P$$

*defined by*

$$\widetilde{\chi_\mathcal{F}^{\text{ind}}}\left([\alpha_n]\right) := \frac{\chi_\mathcal{F}(\alpha_n)}{p_{01} \cdot p_{12} \cdot p_{23} \cdots p_{(n-1)n}}.$$

*Here* $p_{01} := 1$ *and* $\mathcal{R}_P$ *is the ring* $\mathcal{R}_S$ *of fractions of* $\mathcal{R}$ *with respect to the multiplicatively closed set* $S$ *consisting of all the finite products of powers of elements in* $P$.

*(iii) In particular, in the case when the above projective system* $P = \{p^s\}$ *consists of powers of a non-zero element* $p$, *we get the homomorphism*

$$\widetilde{\chi_\mathcal{F}^{\text{ind}}} : \mathcal{F}_P^{\text{st.ind}}(X_\infty) \to \mathcal{R}\left[\frac{1}{p}\right]$$

*defined by*

$$\widetilde{\chi_\mathcal{F}^{\text{ind}}}\left([\alpha_n]\right) := \frac{\chi_\mathcal{F}(\alpha_n)}{p^{n-1}}.$$

*Here* $\mathcal{R}\left[\frac{1}{p}\right]$ *is the ring of polynomials generated by* $\{\frac{1}{p^s}\}$.

*Proof.* (i) follows from taking the inductive limit of the commutative diagram

$$\begin{array}{ccc} \mathcal{F}_P^{\text{st}}(X_\lambda) & \xrightarrow{\chi_\mathcal{F}} & \mathcal{R} \\ \pi_{\lambda\mu}{}^* \downarrow & & \downarrow \times p_{\lambda\mu} \\ \mathcal{F}_P^{\text{st}}(X_\mu) & \xrightarrow{\chi_\mathcal{F}} & \mathcal{R}. \end{array}$$

For a general directed set $\Lambda$, we do not know how to describe the homomorphism $\chi_\mathcal{F}^{\text{ind}}$ in a bit more down-to-earth way. However, when it comes to the case when $\Lambda = \mathbb{N}$, we can get the above claim as follows.

(ii) Let $R_n = \mathcal{R}$ for each $n$ and for $n < m$ let $\rho_{nm} : R_n \to R_m$ denote the homomorphism defined by $\rho_{nm}(r_n) = r_n \cdot p_{n(n+1)} \cdot p_{(n+1)(n+2)} \cdots p_{(m-1)m}$. And let $\phi^n : \mathcal{R}_n \to \mathcal{R}_P$ be the homomorphism defined by

$$\phi^n(r_n) := \frac{r_n}{p_{01} \cdot p_{12} \cdot p_{23} \cdot \cdots \cdot p_{(n-1)n}}.$$

Then we have that for $n < m$

$$\phi^m \circ \rho_{nm} = \phi^n.$$

Therefore it follows from the standard facts of the inductive limits that there exists a unique homomorphism $\Phi : \varinjlim_n \mathcal{R}_n \to \mathcal{R}_P$ such that the following diagram commutes:

$$\begin{array}{ccc} & R_n & \\ {}^{\rho^n}\swarrow & & \searrow{}^{\phi_n} \\ \varinjlim_n R_n & \xrightarrow{\Phi} & \mathcal{R}_P. \end{array}$$

This homomorphism $\Phi : \varinjlim_n R_n \to \mathcal{R}_P$ is a kind of "realization homomorphism" of the abstract ring $\varinjlim_n R_n$. By composing $\chi_{\mathcal{F}}^{\text{ind}} : \mathcal{F}_P^{\text{st.ind}}(X_\infty) \to \varinjlim_n \left\{ \times p_{nm} : \mathcal{R} \to \mathcal{R} \right\}$ with this "realization homomorphism" $\Phi$, we get the above homomorphism $\widetilde{\chi_{\mathcal{F}}^{\text{ind}}} : \mathcal{F}_P^{\text{st.ind}}(X_\infty) \to \mathcal{R}_P$. $\square$

**Remark 3.5.** (i) Let $X_\lambda = pt$ be a point for any $\lambda \in \Lambda$ and let $\pi_{\lambda\mu} = \text{id} : X_\lambda \to X_\mu$ be the identity. Then the proalgebraic variety $\varprojlim\{X_\lambda, \pi_{\lambda\mu} : X_\mu \to X_\lambda\}$ is a point and is called a *proalgebraic point* and denoted by $pt_\infty$. Then for the proalgebraic point $pt_\infty$

$$\mathcal{F}_P^{\text{st.ind}}(pt_\infty) := \varinjlim_{\lambda \in \Lambda} \Big\{ \times p_{\lambda\mu} : \mathcal{F}(pt) \to \mathcal{F}(pt) \Big\}$$
$$= \varinjlim_{\lambda \in \Lambda} \Big\{ \times p_{\lambda\mu} : \mathcal{R} \to \mathcal{R} \Big\}$$

In this sense, the above proalgebraic $\mathcal{F}$-characteristic homomorphism $\chi_{\mathcal{F}}^{\text{ind}} : \mathcal{F}_P^{\text{st.ind}}(X_\infty) \to \varinjlim_{\lambda \in \Lambda} \Big\{ \times p_{\lambda\mu} : \mathcal{R} \to \mathcal{R} \Big\}$ is expressed as $\chi_{\mathcal{F}}^{\text{ind}} : \mathcal{F}_P^{\text{st.ind}}(X_\infty) \to \mathcal{F}_P^{\text{st.ind}}(pt_\infty)$ and it is a proalgebraic version of the $\mathcal{F}$-characteristic $\chi_{\mathcal{F}} : \mathcal{F}(X) \to \mathcal{F}(pt) = \mathcal{R}$.

(ii) Note that if at least one $p_{nm} = 0$ in (ii) above, then $\varinjlim_n \Big\{ \times p_{nm} : \mathcal{R} \to \mathcal{R} \Big\} = 0$; so we assume that all $p_{nm} \neq 0$ in the above theorem.

(iii) The above realization is a *canonical* one in the sense that there are many other realizations by considering other $\phi'_n(r_n) = \dfrac{r_n}{\omega \cdot p_{01} \cdot p_{12} \cdot p_{23} \cdots p_{(n-1)n}}$ with any non-zero element $\omega$.

Certainly here we should discuss the functoriality of $\mathcal{F}_P^{\text{st.ind}}$ for proalgebraic varieties, but we postpone it to §5. Instead, in this section we just discuss the proalgebraic $\chi_{\mathcal{F}}$-characteristic.

In the case when $\mathcal{F} = F$ is the constructible function functor, in the above Theorem 3.4 the ring $\mathcal{R}$ is simply replaced by the integer ring $\mathbb{Z}$. We give some examples:

**Example 3.6.** Let us consider the infinite countable product $X_\infty := X^{\mathbb{N}}$ of a complex algebraic variety $X$ as a simple model case. Let $X^n$ denote the Cartesian product of $n$ copies of the variety $X$. For each projection $\pi_{n(n+1)} : X^{n+1} \to X^n$ (projecting to the first $n$ factors), the pullback homomorphism $\pi_{n(n+1)}^* : F(X^n) \to F(X^{n+1})$ is the multiplication by the characteristic function $\mathbb{1}_X$ of the last factor $X$, i.e.,

$$\pi_{n(n+1)}^*(\alpha) = \alpha \times \mathbb{1}_X,$$

where $(\alpha \times \mathbb{1}_X)(y,x) := \alpha(y)\mathbb{1}_X(x) = \alpha(y)$. Since $\chi(\alpha \times \beta) = \chi(\alpha)\chi(\beta)$, we get the commutative diagram

$$\begin{CD} F(X^n) @>\chi>> \mathbb{Z} \\ @V{\times \mathbb{1}_X}VV @VV{\times \chi(X)}V \\ F(X^{n+1}) @>\chi>> \mathbb{Z}. \end{CD}$$

We assume that $\chi(X) \neq 0$ and let $p_{nm} := \chi(X)^{m-n}$ for $n < m$. Then $P := \{p_{nm}\}$ is a projective system of integers and $F_P^{\text{st.ind}}(X_\infty) = F^{\text{ind}}(X_\infty)$ and we get the proalgebraic Euler–Poincaré characteristic $\chi^{\text{ind}} : F^{\text{ind}}(X_\infty) \to \mathbb{Z}\left[\frac{1}{\chi(X)}\right]$ described by

$$\chi^{\text{ind}}([\alpha_n]) = \frac{\chi(\alpha_n)}{\chi(X)^{n-1}}.$$

In particular, we have that $\chi^{\text{ind}}([\mathbb{1}_{X^n}]) = \chi(X)$ for any $n$, since $\chi(X^n) = \chi(X)^n$.

**Example 3.7.** Let $X_\infty = \varprojlim_{n \in \mathbb{N}}\left\{X_n, \pi_{nm} : X_m \to X_n\right\}$ be a proalgebraic variety such that for each $n$ the structure morphism $\pi_{n(n+1)} : X_{n+1} \to X_n$ satisfies the condition that the Euler–Poincaré characteristics of the fibers of $\pi_{n(n+1)}$ are non-zero (which implies the surjectivity of the morphism $\pi_{n(n+1)}$) and constant; for example, $\pi_{n(n+1)} : X_{n+1} \to X_n$ is a locally trivial fiber bundle with fiber variety being $F_n$ and $\chi(F_n) \neq 0$. Let us denote the constant Euler–Poincaré characteristic of the fibers of the morphism $\pi_{n(n+1)} : X_{n+1} \to X_n$ by $e_n$ and we set $e_0 := 1$. Then we get the canonical proalgebraic Euler–Poincaré characteristic homomorphism

$$\chi^{\text{ind}} : F^{\text{ind}}(X_\infty) \to \mathbb{Q}$$

described by

$$\chi^{\text{ind}}([\alpha_n]) = \frac{\chi(\alpha_n)}{e_0 \cdot e_1 \cdot e_2 \cdots e_{n-1}}.$$

In particular, if the Euler–Poincaré characteristics $e_n$ are all the same, say $e_n = e$ for any $n$, then the canonical proalgebraic Euler–Poincaré characteristic homomorphism $\chi^{\text{ind}} : F^{\text{ind}}(X_\infty) \to \mathbb{Q}$ is described by $\chi^{\text{ind}}([\alpha_n]) = \dfrac{\chi(\alpha_n)}{e^{n-1}}$, and furthermore the target ring $\mathbb{Q}$ can be replaced by the ring $\mathbb{Z}\left[\frac{1}{e}\right]$.

In this example, we need the commutativity of $\chi$ with the pushforward, although we do not use the commutativity of $\chi_\mathcal{F}$ with the pushforward in the above Theorem 3.4. Let $f : X \to Y$ be a morphism such that its fibers all have the same non-zero Euler–Poincaré characteristic, denoted by $e_f$. Then we can see that for any characteristic function $\mathbb{1}_W$ we have

$$f_* f^* \mathbb{1}_W = e_f \cdot \mathbb{1}_W.$$

Hence we have

$$\begin{aligned} \chi(f^*\mathbb{1}_W) &= \chi(f_*f^*\mathbb{1}_W) \\ &= \chi(e_f \cdot \mathbb{1}_W) \\ &= e_f \cdot \chi(\mathbb{1}_W). \end{aligned}$$

Therefore for any constructible function $\alpha \in F(Y)$ we get $\chi(f^*\alpha) = e_f \cdot \chi(\alpha)$. Hence if we set

$$p_{nm} = \begin{cases} 1 & n = m \\ e_n \cdot e_{n+1} \cdots e_{m-1} & n < m, \end{cases}$$

then $P := \{p_{nm}\}$ is a projective system and $F_P^{\text{st.ind}}(X_\infty) = F^{\text{ind}}(X_\infty)$. Thus the above description of $\chi^{\text{ind}}([\alpha_n])$ follows from the above theorem.

Example 3.7 motivates us to define the following notion.

**Definition 3.8.** Let $\mathcal{F}$ be a bifunctor on a category $\mathcal{C}$ such that $\mathcal{R} = \mathcal{F}(pt)$ is a commutative ring with a unit and let $\chi_{\mathcal{F}} : \mathcal{F}(X) \to \mathcal{R}$ be the $\mathcal{F}$-characteristic.
(3.8.1) If a morphism $f : X \to Y$ satisfies the condition that for an element $\alpha \in \mathcal{F}(Y)$
$$\chi_{\mathcal{F}}(f_* f^* \alpha) = c_f \cdot \chi_{\mathcal{F}}(\alpha)$$
with some multiplier $c_f \in \mathcal{R}$ depending only on the morphism $f$, then we say that $f$ is $\chi_{\mathcal{F}}$-*constant with respect to* $\alpha$ *with the multipler* $c_f$. ($c_f$ could be considered as the "$\chi_F$-characteristic of the fiber of $f$".)
(3.8.2) If $f$ is a $\chi_{\mathcal{F}}$-*constant with respect to any element* $\alpha \in \mathcal{F}(Y)$ *with the multipler* $c_f$, then the morphism $f : X \to Y$ is called $\chi_{\mathcal{F}}$-*constant with the multipler* $c_f$.
(3.8.3) (a bit stronger than (3.8.2)) Let $\mathcal{F}$ be a bifunctor from a category $\mathcal{C}$ to the category of $\mathcal{R}$-modules such that $\mathcal{F}(pt) = \mathcal{R}$. If a morphism $f : X \to Y$ satisfies the condition that
$$f_* f^* = c_f \cdot Id_{\mathcal{F}(Y)} : \mathcal{F}(Y) \to \mathcal{F}(Y)$$
with some element $c_f \in \mathcal{R}$, where $Id_{\mathcal{F}(Y)}$ denotes the identity homomorphism, then $f$ is also called $\chi_{\mathcal{F}}$-constant with the multiplier $c_f$. (Note that in this case $f_* f^* = c_f \cdot Id_{\mathcal{F}(Y)}$ implies that $\chi_{\mathcal{F}}(f_* f^* \alpha) = c_f \cdot \chi_{\mathcal{F}}(\alpha)$ for any $\alpha \in \mathcal{F}(Y)$.)

Examples of a bifunctor $\mathcal{F}$ from a category $\mathcal{C}$ to the category of $\mathcal{R}$-modules such that $\mathcal{F}(pt) = \mathcal{R}$ are *Green functors*, which are discussed later in §7.

So, with this definition, $\alpha_\lambda \in \mathcal{F}_P^{\text{st}}(X_\lambda)$ means that $\pi_{\lambda\mu}$ is $\chi_{\mathcal{F}}$-constant with respect to $\alpha_\lambda$ with the multiplier $p_{\lambda\mu}$ for any $\mu$ such that $\lambda < \mu$.

With this definition, the above Example 3.7 can be generalized to the following

**Example 3.9.** Let $X_\infty = \varprojlim_{n \in \mathbb{N}} \left\{ X_n, \pi_{nm} : X_m \to X_n \right\}$ be a proalgebraic variety such that for each $n$ the structure morphism $\pi_{n(n+1)} : X_{n+1} \to X_n$ is $\chi_{\mathcal{F}}$-constant with the multiplier $c_{n(n+1)} \in \mathcal{R}$. Then we get the canonical proalgebraic $\mathcal{F}$-characteristic homomorphism $\chi_{\mathcal{F}}^{\text{ind}} : \mathcal{F}^{\text{ind}}(X_\infty) \to \mathcal{R}_P$ defined by
$$\chi_{\mathcal{F}}^{\text{ind}}\left( [\alpha_n] \right) := \frac{\chi_{\mathcal{F}}(\alpha_n)}{c_{01} \cdot c_{12} \cdot c_{23} \cdots c_{(n-1)n}}.$$
Here $c_{01} := 1$ and $\mathcal{R}_P$ is the ring $\mathcal{R}_S$ of fractions of $\mathcal{R}$ with respect to the multiplicatively closed set $S$ consisting of all the finite products of powers of multipliers $\{c_{nm}\}$.

We can show other examples, using Fulton–MacPherson's Bivariant Theory [MF] (also see [F1]). So, we quickly recall only necessary ingredients of the Bivariant Theory for using it in this section and later sections.

A bivariant theory $\mathbb{B}$ on a category $\mathcal{C}$ with values in the category of abelian groups is an assignment to each morphism
$$X \xrightarrow{f} Y$$
in the category $\mathcal{C}$ an abelian group
$$\mathbb{B}(X \xrightarrow{f} Y)$$
which is equipped with the following three basic operations:
**Products**: For morphisms $f : X \to Y$ and $g : Y \to Z$, the product operation
$$\bullet : \mathbb{B}(X \xrightarrow{f} Y) \otimes \mathbb{B}(Y \xrightarrow{g} Z) \to \mathbb{B}(X \xrightarrow{gf} Z),$$
**Pushforwards**: For morphisms $f : X \to Y$ and $g : Y \to Z$ with $f$ proper, the pushforward operation
$$f_\star : \mathbb{B}(X \xrightarrow{gf} Z) \to \mathbb{B}(Y \xrightarrow{g} Z)$$

**Pullbacks**: For a fiber square

$$\begin{array}{ccc} X' & \xrightarrow{g'} & X \\ f' \downarrow & & \downarrow f \\ Y' & \xrightarrow{g} & Y, \end{array}$$

the pullback operation

$$g^\star : \mathbb{B}(X \xrightarrow{f} Y) \to \mathbb{B}(X' \xrightarrow{f'} Y').$$

And these three operations are required to satisfy the seven compatibility axioms (see [FM, Part I, §2.2] for details).

A bivariant theory $\mathbb{B}$ is said to *have units* (see [FM, §2.2]) if there exists an element $1_X \in \mathbb{B}(X \xrightarrow{id_X} X)$ such that $\alpha \bullet 1_X = \alpha$ for all maps $W \to X$ and all $\alpha \in \mathbb{B}(W \to X)$; that $1_X \bullet \beta = \beta$ for all maps $X \to Y$ and all $\beta \in \mathbb{B}(X \to Y)$; and that $g^\star 1_X = 1_{X'}$ for all $g : X' \to X$.

Let $\mathbb{B}, \mathbb{B}'$ be two bivariant theories on a category $\mathcal{C}$. Then a *Grothendieck transformation* from $\mathbb{B}$ to $\mathbb{B}'$

$$\gamma : \mathbb{B} \to \mathbb{B}'$$

is a collection of homomorphisms

$$\mathbb{B}(X \to Y) \to \mathbb{B}'(X \to Y)$$

for a morphism $X \to Y$ in the category $\mathcal{C}$, which preserve the above three basic operations:
   (i)   $\gamma(\alpha \bullet_\mathbb{B} \beta) = \gamma(\alpha) \bullet_{\mathbb{B}'} \gamma(\beta)$,
   (ii)  $\gamma(f_\star \alpha) = f_\star \gamma(\alpha)$, and
   (iii) $\gamma(g^\star \alpha) = g^\star \gamma(\alpha)$.

A bivariant theory unifies both a covariant theory and a contravariant theory in the following sense: $\mathbb{B}_*(X) := \mathbb{B}(X \to pt)$ and $\mathbb{B}^*(X) := \mathbb{B}(X \xrightarrow{id} X)$ become a covariant functor and a contravariant functor, respectively. And a Grothendieck transformation $\gamma : \mathbb{B} \to \mathbb{B}'$ induces natural transformations $\gamma_* : \mathbb{B}_* \to \mathbb{B}'_*$ and $\gamma^* : \mathbb{B}^* \to \mathbb{B}'^*$.

For the sake of a later use, we also note that if we have a Grothendieck transformation $\gamma : \mathbb{B} \to \mathbb{B}'$, then via a bivariant class $b \in \mathbb{B}(X \xrightarrow{f} Y)$ we get the commutative diagram

$$\begin{array}{ccc} \mathbb{B}_*(Y) & \xrightarrow{\gamma_*} & \mathbb{B}'_*(Y) \\ b\bullet \downarrow & & \downarrow \gamma(b)\bullet \\ \mathbb{B}_*(X) & \xrightarrow{\gamma_*} & \mathbb{B}'_*(X). \end{array}$$

This is called *the Verdier-type Riemann–Roch formula associated to the bivariant class $b$*.

Fulton–MacPherson's bivariant group $\mathbb{F}(X \xrightarrow{f} Y)$ of constructible functions consists of all the constructible functions on $X$ which satisfy the local Euler condition with respect to $f$. Here a constructible function $\alpha \in F(X)$ is said to satisfy the *local Euler condition with respect to $f$* if for any point $x \in X$ and for any local embedding $(X, x) \to (\mathbf{C}^N, 0)$ the equality $\alpha(x) = \chi\left(B_\epsilon \cap f^{-1}(z); \alpha\right)$ holds, where $B_\epsilon$ is a sufficiently small open ball of the origin $0$ with radius $\epsilon$ and $z$ is any point close to $f(x)$ (cf. [Br], [Sa2]). In particular, if $\mathbb{1}_f := \mathbb{1}_X$ belongs to the bivariant group $\mathbb{F}(X \xrightarrow{f} Y)$, then the morphism $f : X \to Y$ is called *an Euler morphism*. And any constructible function in the bivariant group $\mathbb{F}(X \xrightarrow{f} Y)$ is called *a bivariant constructible function* to emphasize the bivariantness.

The three operations on $\mathbb{F}$ are defined as follows:
(i) the product $\bullet : \mathbb{F}(X \xrightarrow{f} Y) \otimes \mathbb{F}(Y \xrightarrow{g} Z) \to \mathbb{F}(X \xrightarrow{gf} Z)$ is defined by

$$\alpha \bullet \beta := \alpha \cdot f^* \beta,$$

(ii) the pushforward $f_\star : \mathbb{F}(X \xrightarrow{gf} Z) \to \mathbb{F}(Y \xrightarrow{g} Z)$ is the usual pushforward $f_*$, i.e.,
$$f_\star(\alpha)(y) := \int c_*(\alpha|_{f^{-1}}),$$

(iii) for a fiber square
$$\begin{array}{ccc} X' & \xrightarrow{g'} & X \\ f' \downarrow & & \downarrow f \\ Y' & \xrightarrow{g} & Y, \end{array}$$

the pullback $g^\star : \mathbb{F}(X \xrightarrow{f} Y) \to \mathbb{F}(X' \xrightarrow{f'} Y')$ is the functional pullback $g'^*$, i.e..,
$$g^\star(\alpha)(x') := \alpha(g'(x')).$$

Note that $\mathbb{F}(X \xrightarrow{\mathrm{id}_X} X)$ consists of all locally constant functions and $\mathbb{F}(X \to pt) = F(X)$. As a corollary of this observation, we have

**Proposition 3.10.** *For any bivariant constructible function $\alpha \in \mathbb{F}(X \xrightarrow{f} Y)$, the Euler–Poincaré characteristic $\chi(f^{-1}(y); \alpha) = \int c_*(\alpha|_{f^{-1}(y)})$ of $\alpha$ restricted to each fiber $f^{-1}(y)$ is locally constant, i.e., constant along connected components of the base variety $Y$. In particular, if $f : X \to Y$ is an Euler proper morphism, then the Euler–Poincaré characteristic of the fibers are locally constant.*

Note that locally trivial fiber bundles are Euler, but not vice versa.

**Example 3.11.** Let $X_\infty = \varprojlim_{\lambda \in \Lambda} \{X_\lambda, \pi_{\lambda\mu} : X_\mu \to X_\lambda\}$ be a proalgebraic variety such that for each $\lambda < \mu$ the structure morphism $\pi_{\lambda\mu} : X_\mu \to X_\lambda$ is an Euler proper morphism (hence surjective) of topologically connected algebraic varieties with the constant Euler–Poincaré characteristic $\chi_{\lambda\mu}$ of the fiber of the morphism $\pi_{\lambda\mu}$ being non-zero. Then we get the proalgebraic Euler–Poincaré characteristic homomorphism
$$\chi^{\mathrm{ind}} : F^{\mathrm{ind}}(X_\infty) \to \varinjlim_\Lambda \{\times \chi_{\lambda\mu} : \mathbb{Z} \to \mathbb{Z}\}.$$

For a morphism $f : X \to Y$ and a bivariant class $b \in \mathbb{B}(X \xrightarrow{f} Y)$, the pair $(f; b)$ is called a *bivariant-class-equipped morphism* and we just express $(f; b) : X \to Y$. Let $\mathbb{B}$ be a bivariant theory having units. If a system $\{b_{\lambda\mu}\}$ of bivariant classes satisfies that
$$b_{\lambda\lambda} = 1_{X_\lambda} \quad \text{and} \quad b_{\mu\nu} \bullet b_{\lambda\mu} = b_{\lambda\nu} \quad (\lambda < \mu < \nu),$$
then we call the system *a projective system of bivariant classes*. If $\{\pi_{\lambda\mu} : X_\mu \to X_\lambda\}$ and $\{b_{\lambda\mu}\}$ are projective systems, then the system $\{(\pi_{\lambda\mu}; b_{\lambda\mu}) : X_\mu \to X_\lambda\}$ shall be called *a projective system of bivariant-class-equipped morphisms*.

For a bivariant theroy $\mathbb{B}$ having units on the category $\mathcal{C}$ and for a projective system $\{(\pi_{\lambda\mu}; b_{\lambda\mu}) : X_\mu \to X_\lambda\}$ of bivariant-class-equipped morphisms, the inductive limit
$$\varinjlim_\Lambda \Big\{ \mathbb{B}_*(X_\lambda), b_{\lambda\mu} \bullet : \mathbb{B}_*(X_\lambda) \to \mathbb{B}_*(X_\mu) \Big\}$$
shall be denoted by
$$\mathbb{B}_*^{\mathrm{ind}}\Big(X_\infty; \{b_{\lambda\mu}\}\Big)$$
emphasizing the projective system $\{b_{\lambda\mu}\}$ of bivariant classes, because the above inductive limit surely depends on the choice of it. For example, in the above Example 3.11 we have that
$$F^{\mathrm{ind}}(X_\infty) = F_*^{\mathrm{ind}}\Big(X_\infty; \{\mathbb{1}_{\pi_{\lambda\mu}}\}\Big).$$

**Example 3.12.** Let $\left\{(\pi_{n(n+1)}, \alpha_{n(n+1)}) : X_{n+1} \to X_n\right\}$ be a projective system of bivariant - class - equipped morphisms of topologically connected algebraic varieties with $\alpha_{n(n+1)} \in \mathbb{F}(X_{n+1} \to X_n)$. And assume that the (constant) Euler–Poincaré characteristic $\chi\left(\pi_{n(n+1)}{}^{-1}(y); \alpha_{n(n+1)}\right)$ of $\alpha_{n(n+1)}$ restricted to each fiber $\pi_{n(n+1)}{}^{-1}(y)$ is non-zero and it shall be denoted by $e_f(\alpha_{n(n+1)})$. And we set $e_f(\alpha_{01}) := 1$. Then the canonical Euler–Poincaré characteristic homomorphism

$$\chi^{\text{ind}} : F_*^{\text{ind}}\left(X_\infty; \{\alpha_{n(n+1)}\}\right) \to \mathbb{Q}$$

is described by

$$\chi^{\text{ind}}([\alpha_n]) = \frac{\chi(\alpha_n)}{e_f(\alpha_{01}) \cdot e_f(\alpha_{12}) \cdots e_f(\alpha_{(n-1)n})}.$$

This can be seen as follows. Let $(f, \alpha) : X \to Y$ be a bivariant-class-equipped morphism of topologically connected algebraic varieties with $\alpha \in \mathbb{F}(X \xrightarrow{f} Y)$. It follows from Proposition 3.10 that the Euler–Poincaré characteristic $\chi\left(f^{-1}(y); \alpha\right)$ of $\alpha$ restricted to each fiber $f^{-1}(y)$ is constant (and non-zero by assumption). So, if it is denoted by $e_f(\alpha)$, then $f_*\alpha = e_f(\alpha) \cdot \mathbb{1}_Y$. Then to prove the above statement, it suffices to see that we have the following commutative diagram:

$$\begin{array}{ccc} F(Y) & \xrightarrow{\chi} & \mathbb{Z} \\ \alpha \bullet \downarrow & & \downarrow \times e_f(\alpha) \\ F(X) & \xrightarrow{\chi} & \mathbb{Z}. \end{array}$$

To see this, we need the *projection formula* that for a morphism $f : X \to Y$ and constructible functions $\alpha \in F(X)$ and $\beta \in F(Y)$

$$f_*(\alpha \cdot f^*\beta) = (f_*\alpha) \cdot \beta.$$

Then, using this projection formula we have

$$\begin{aligned} \chi(\alpha \bullet \beta) &= \chi(\alpha \cdot f^*\beta) \\ &= \chi(f_*\alpha \cdot \beta) \\ &= \chi\left((e_f(\alpha) \cdot \mathbb{1}_Y) \cdot \beta\right) \\ &= e_f(\alpha) \cdot \chi(\beta) \end{aligned}$$

Thus we get the above commutative diagram.

To get a similar result in the above more general case of $\mathbb{B}_*^{\text{ind}}\left(X_\infty; \{b_{\lambda\mu}\}\right)$ we assume that $\mathbb{B}_*(pt)$ is a commutative ring with a unit, denoted by $\mathcal{R}^\mathbb{B}$, and let $P = \{p_{\lambda\mu}\}$ be a projective system of non-zero elements $p_{\lambda\mu} \in \mathcal{R}^\mathbb{B}$. Then, if we set

$$\mathbb{B}_{*,P}^{\text{st.ind}}\left(X_\infty; \{b_{\lambda\mu}\}\right)$$
$$:= \left\{[\alpha_\lambda] \in \mathbb{B}_*^{\text{ind}}\left(X_\infty; \{b_{\lambda\mu}\}\right) | \chi_{\mathbb{B}_*}(b_{\lambda\mu} \bullet \alpha_\lambda) = p_{\lambda\mu} \cdot \chi_{\mathbb{B}_*}(\alpha_\lambda) \quad (\lambda < \mu)\right\},$$

we get the proalgebraic $\chi_{\mathbb{B}_*}$-characteristic homomorphism

$$\chi_{\mathbb{B}_*}^{\text{ind}} : \mathbb{B}_{*,P}^{\text{st.ind}}\left(X_\infty; \{b_{\lambda\mu}\}\right) \to \varinjlim_{\lambda \in \Lambda}\left\{\times p_{\lambda\mu} : \mathcal{R}^\mathbb{B} \to \mathcal{R}^\mathbb{B}\right\}.$$

4. Relative Grothendieck rings and motivic measures

In the previous section we have dealt with the constructible function functor $F$ as an example of a bifunctor $\mathcal{F}$. In this section, we will deal with the so-called *relative Grothendieck ring of complex algebraic varieties over $X$*, denoted by $K_0(\mathcal{V}/X)$. This was introduced by E. Looijenga in [Lo] and further studied by F. Bittner in [Bi]. For a very recent application of the relative Grothendieck groups, see [BSY2] (cf. [BSY2']).

The relative Grothendieck group $K_0(\mathcal{V}/X)$ is the quotient of the free abelian group of isomorphism classes of morphisms to $X$ (denoted by $[Y \to X]$ or $[Y \xrightarrow{h} X]$), modulo the following relation:

$$[Y \xrightarrow{h} X] = [Z \hookrightarrow Y \xrightarrow{h} X] + [Y \setminus Z \hookrightarrow Y \xrightarrow{h} X]$$

for $Z \subset Y$ a closed subvariety of $Y$. The ring structure is given by the fiber square: for $[Y \xrightarrow{f} X], [W \xrightarrow{g} X] \in K_0(\mathcal{V}/X)$

$$[Y \xrightarrow{f} X] \cdot [W \xrightarrow{g} X] := [Y \times_X W \xrightarrow{f \times_X g} X].$$

Here $Y \times_X W \xrightarrow{f \times_X g} X$ is $g \circ f' = f \circ g'$ where $f'$ and $g'$ are as in the following diagram

$$\begin{array}{ccc} Y \times_X W & \xrightarrow{f'} & W' \\ g' \downarrow & & \downarrow g \\ Y & \xrightarrow{f} & X. \end{array}$$

The relative Grothendieck ring $K_0(\mathcal{V}/X)$ has the unit $1_X := [X \xrightarrow{id_X} X]$.

When $X = pt$ is a point, the relative Grothendieck ring $K_0(\mathcal{V}/pt)$ is the usual Grothendieck ring $K_0(\mathcal{V})$ of $\mathcal{V}$, i.e., the free abelian group generated by the isomorphism classes of varieties modulo the subgroup generated by elements of the form $[V] - [V'] - [V \setminus V']$ for a subvariety $V' \subset V$, and the ring structure is given by the Cartesian product of varieties.

For a morphism $f: X' \to X$, the pushforward

$$f_*: K_0(\mathcal{V}/X') \to K_0(\mathcal{V}/X)$$

is defined by

$$f_*[Y \xrightarrow{h} X'] := [Y \xrightarrow{f \circ h} X].$$

With this pushforward, the assignment $X \longmapsto K_0(\mathcal{V}/X)$ is a covariant functor. The pullback

$$f^*: K_0(\mathcal{V}/X) \to K_0(\mathcal{V}/X')$$

is defined as follows: for a fiber square

$$\begin{array}{ccc} Y' & \xrightarrow{g'} & X' \\ f' \downarrow & & \downarrow f \\ Y & \xrightarrow{g} & X \end{array}$$

the pullback $f^*[Y \xrightarrow{g} X] := [Y' \xrightarrow{g'} X']$. With this pullback, the assignment $X \longmapsto K_0(\mathcal{V}/X)$ is a contravariant functor.

Hence, the assignment $X \longmapsto K_0(\mathcal{V}/X)$ is a bifunctor and just like in the constructible function functor $F$, by considering the map to a point $\pi_X : X \to pt$, we get the following homomorphism

$$\pi_{X*}: K_0(\mathcal{V}/X) \to K_0(\mathcal{V})$$

which shall be denoted by $\chi_{\mathrm{Gro}}$. And also we get the following commutative diagram:

$$\begin{array}{ccc} K_0(\mathcal{V}/X) & & \\ {\scriptstyle f_*}\downarrow & \searrow^{\chi_{\mathrm{Gro}}} & \\ & & K_0(\mathcal{V}) \\ & \nearrow_{\chi_{\mathrm{Gro}}} & \\ K_0(\mathcal{V}/Y) & & \end{array}$$

and we get the same results as in Theorem 3.4 by replacing $\mathcal{F}(?)$ by $K_0(\mathcal{V}/?)$.

**Observation 4.1.** *A Zariski locally trivial fiber bundle is a $\chi_{\mathrm{Gro}}$-constant morphism with the multiplier being the Grothendieck class of its fiber variety.*

There exists a canonical homomorphism $e : K_0(\mathcal{V}/X) \to F(X)$ (see [BSY2]) defined by
$$e([Y \xrightarrow{f} X]) := f_* \mathbb{1}_Y,$$
which is compatible with the pushforward, i.e., the correspondence of covariant functors $e : K_0(\mathcal{V}/?) \to F(?)$ is a natural transformation. It will be explained in §6 that this natural transformation is unique in a sense.

There exists a canonical homomorphism $\iota : F(X) \to K_0(\mathcal{V}/X)$ defined by $\iota(\mathbb{1}_W) := [W \xrightarrow{i_W} X]$, where $i_W : W \to X$ is the inclusion map. The composite homomorphism $\Gamma := \chi_{\mathrm{Gro}} \circ \iota : F(X) \to K_0(\mathcal{V})$ is more directly and simply defined by
$$\Gamma(\mathbb{1}_W) := [W] \quad \text{or more meaningfully} \quad \Gamma(\alpha) = \sum_{n \in \mathbb{Z}} n \left[ \alpha^{-1}(n) \right].$$

And we have the following commutative diagram:

$$\begin{array}{ccc} F(X) & \xrightarrow{\Gamma} & K_0(\mathcal{V}) \\ & \searrow_{\chi} \quad \swarrow_{e} & \\ & \mathbb{Z} & \end{array}$$

**Definition 4.2.** Let $R$ be a commutative ring. A map
$$E : \mathrm{Obj}(\mathcal{V}) \to R$$
is called a *generalized Euler characteristic with value in $R$* if the following three conditions hold:
 (i) $E(X) = E(Y)$ if $X \cong Y$,
 (ii) $E(X) = E(Y) + E(X \setminus Y)$ for $Y \subset X$,
 (iii) $E(X \times Y) = E(X) \cdot E(Y)$.

A typical example of $E$ is of course the topological Euler–Poincaré characteristic $\chi$ with $R = \mathbb{Z}$ and $E$ induces the homomorphism $E_F : F(X) \to R$ defined simply by $E_F(\sum_S a_S \mathbb{1}_S) := \sum_S a_S E(S)$. And $E_F$ factors through the above "tautological" homomorphism $\Gamma : F(X) \to K_0(\mathcal{V})$:

$$\begin{array}{ccc} F(X) & \xrightarrow{\Gamma} & K_0(\mathcal{V}) \\ & \searrow_{E_F} \quad \swarrow_{\widetilde{E}} & \\ & R & \end{array}$$

where $\widetilde{E} : K_0(\mathcal{V}) \to \mathcal{R}$ is defined by $\widetilde{E}([X]) := E(X)$.

So $\Gamma : F(X) \to K_0(\mathcal{V})$ is a "motivic" version of the topological Euler–Poincaré characteristic $\chi : F(X) \to \mathbb{Z}$ and provisionally called the *Grothendieck class homomorphism*.

In the previous section we have generalized $\chi : F(X) \to \mathbb{Z}$ to the category of proalgebraic varieties. It turns out that in a similar way as in the previous section we can generalize the Grothendieck class homomorphism $\Gamma : F(X) \to K_0(\mathcal{V})$ to the category of proalgebraic varieties. Here we emphasize that unlike the Euler–Poincaré characteristic $\chi$, $\Gamma : F(X) \to K_0(\mathcal{V})$ is not compatible with the pushforward $f_* : F(X) \to F(Y)$ for a morphism $f : X \to Y$, i.e., the following diagram is *not* commutative:

$$\begin{array}{c} F(X) \\ f_* \downarrow \quad \searrow^{\Gamma} \\ \quad \quad \quad K_0(\mathcal{V}). \\ \quad \quad \nearrow_{\Gamma} \\ F(Y) \end{array}$$

Let $G = \{\gamma_{\lambda\mu}\}$ be a projective system of non-zero Grothendieck classes $\gamma_{\lambda\mu} \in K_0(\mathcal{V})$ indexed by the directed set $\Lambda$, as in §3. Then in the same way as done in §3, we can define

$$F_G^{\mathrm{st}}(X_\lambda) := \Big\{\alpha_\lambda \in F(X_\lambda) \quad | \quad \Gamma\big(\pi_{\lambda\mu}{}^*\alpha_\lambda\big) = \gamma_{\lambda\mu} \cdot \Gamma(\alpha_\lambda) \quad \text{for any} \quad \mu > \lambda\Big\}.$$

For each $\lambda \in \Lambda$, an element of $F_G^{\mathrm{st}}(X_\lambda)$ is called *a $\Gamma$-stable constructible function with respect to the projective system $G$ of non-zero Grothendieck classes*. And for a proalgebraic variety $X_\infty = \varprojlim_{\lambda \in \Lambda} \Big\{X_\lambda, \pi_{\lambda\mu} : X_\mu \to X_\lambda\Big\}$ we define

$$F_G^{\mathrm{st.ind}}(X_\infty) := \varinjlim_{\lambda \in \Lambda}\Big\{F_G^{\mathrm{st}}(X_\lambda), \quad \pi_{\lambda\mu}{}^* : F_G^{\mathrm{st}}(X_\lambda) \to F_G^{\mathrm{st}}(X_\mu) \quad (\lambda < \mu)\Big\}$$

and an element of this group shall be called a *$\Gamma$-stable indconstructible function on the proalgebraic variety $X_\infty$ with respect to the projective system $G$ of non-zero Grothendieck classes*. And we get the following, which is stated as a theorem:

**Theorem 4.3.** *(i) For a proalgebraic variety $X_\infty = \varprojlim_{\lambda \in \Lambda}\Big\{X_\lambda, \pi_{\lambda\mu} : X_\mu \to X_\lambda\Big\}$ and a projective system $G = \{\gamma_{\lambda\mu}\}$ of non-zero Grothendieck classes, we get the proalgebraic Grothendieck class homomorphism*

$$\Gamma^{\mathrm{ind}} : F_G^{\mathrm{st.ind}}(X_\infty) \to \varinjlim_{\lambda \in \Lambda}\Big\{\times \gamma_{\lambda\mu} : K_0(\mathcal{V}) \to K_0(\mathcal{V})\Big\}.$$

*(ii) In the case when $\Lambda = \mathbb{N}$, for a proalgebraic variety $X_\infty = \varprojlim_{n \in \mathbb{N}}\Big\{X_n, \pi_{nm} : X_m \to X_n\Big\}$ and a projective system $G = \{\gamma_{n,m}\}$ of non-zero Grothendieck classes, we have the following canonical proalgebraic Grothendieck class homomorphism*

$$\widetilde{\Gamma^{\mathrm{ind}}} : F_G^{\mathrm{st.ind}}(X_\infty) \to K_0(\mathcal{V})_G$$

*which is defined by*

$$\widetilde{\Gamma^{\mathrm{ind}}}\Big([\alpha_n]\Big) := \frac{\Gamma(\alpha_n)}{\gamma_{01} \cdot \gamma_{12} \cdot \gamma_{23} \cdots \gamma_{(n-1)n}}.$$

*Here we set $\gamma_{01} := \mathbb{1}$ and $K_0(\mathcal{V})_G$ is the ring of fractions of $K_0(\mathcal{V})$ with respect to the multiplicatively closed set consisting of finite products of powers of elements of $G$.*

*(iii) Let $X_\infty = \varprojlim_{n \in \mathbb{N}}\Big\{X_n, \pi_{nm} : X_m \to X_n\Big\}$ be a proalgebraic variety such that each structure morphism $\pi_{n(n+1)} : X_{n+1} \to X_n$ satisfies the condition that for each $n$ there exists a $\gamma_n \in K_0(\mathcal{V})$ such that $\pi_{n(n+1)}{}^{-1}(S_n) = \gamma_n \cdot [S_n]$ for any constructible set $S_n \subset X_n$; for example, $\pi_{n(n+1)} : X_{n+1} \to X_n$ is a Zariski locally trivial fiber bundle* *with fiber variety being $F_n$ (in which case $\gamma_n = [F_n] \in K_0(\mathcal{V})$). Then the canonical proalgebraic Grothendieck class homomorphism*

$$\Gamma^{\mathrm{ind}} : F^{\mathrm{ind}}(X_\infty) \to K_0(\mathcal{V})_G$$

*is described by*

$$\Gamma^{\mathrm{ind}}([\alpha_n]) = \frac{\Gamma(\alpha_n)}{\gamma_0 \cdot \gamma_1 \cdot \gamma_2 \cdots \gamma_{n-1}}.$$

*Here $\gamma_0 := \mathbb{1}$ and $K_0(\mathcal{V})_G$ is the ring of fractions of $K_0(\mathcal{V})$ with respect to the multiplicatively closed set consisting of finite products of powers of $\gamma_m$ ($m = 1, 2, 3 \cdots$).*

*(iv) In particular, if $\gamma_n$ are all the same, say $\gamma_n = \gamma$ for any $n$, then the canonical proalgebraic Grothendieck class homomorphism*

$$\Gamma^{\mathrm{ind}} : F^{\mathrm{ind}}(X_\infty) \to K_0(\mathcal{V})_G$$

*is described by*

$$\Gamma^{\mathrm{ind}}([\alpha_n]) = \frac{\Gamma(\alpha_n)}{\gamma^{n-1}}.$$

*In this special case the quotient ring $K_0(\mathcal{V})_G$ shall be simply denoted by $K_0(\mathcal{V})_\gamma$.*

**Remark 4.4.** When we consider a localization or a ring of fractions of the Grothendieck ring $K_0(\mathcal{V})$, we need to be a bit careful. Unlike the ring $\mathbb{Z}$ of integers, the Grothendieck ring $K_0(\mathcal{V})$ is not a domain, which is a recent result due to B. Poonen [Po, Theorem 1]. Also it is in general hard to check whether the class Grothendieck class $[V]$ of a variety $V$ is a non-zero divisor or not; indeed, one does not know whether even the Grothendieck class $[\mathbb{P}^n]$ of the projective space is a non-zero divisor (which Willem Veys pointed out to the author).

**Example 4.5.** The arc space $\mathcal{L}(X)$ of an algebraic variety $X$ is defined to be the projective limit of the projective system consisting of truncated arc varieties $\mathcal{L}_n(X)$ and projections $\pi_{n(n+1)} : \mathcal{L}_{n+1}(X) \to \mathcal{L}_n(X)$. Thus the arc space is a nontrivial example of a proalgebraic variety. If $X$ is nonsingular and of complex dimension $d$, then the projection $\pi_{n(n+1)} : \mathcal{L}_{n+1}(X) \to \mathcal{L}_n(X)$ is a Zariski locally trivial fiber bundle with fiber being $\mathbb{C}^d$. Thus in this case, in Theorem 4.3 (iv) the Grothendieck class $\gamma$ is $\mathbb{L}^d$.

The indconstructible function is just an element of $F^{\mathrm{ind}}(X_\infty) = \varinjlim_{\lambda \in \Lambda} F(X_\lambda)$ and up to now we do not discuss the role of function, even though it is called "function". In fact, the indconstructible function can be considered in a natural way as a function on the proalgebraic variety simply as follows: for $[\alpha_\lambda] \in F^{\mathrm{ind}}(X_\infty) = \varinjlim_{\lambda \in \Lambda} F(X_\lambda)$ the value of $[\alpha_\lambda]$ at a point $(x_\mu) \in X_\infty = \varprojlim_{\lambda \in \Lambda} X_\lambda$ is defined by

$$[\alpha_\lambda]\bigl((x_\mu)\bigr) := \alpha_\lambda(x_\lambda)$$

which is well-defined. So, if we let $Fun(X_\infty, \mathbb{Z})$ be the abelian group of $\mathbb{Z}$-valued functions on $X_\infty$, then the homomorphism

$$\Psi : \varinjlim_{\lambda \in \Lambda} F(X_\lambda) \to Fun(X_\infty, \mathbb{Z}) \quad \text{defined by} \quad \Psi([\alpha_\lambda])((x_\mu)) := \alpha_\lambda(x_\lambda)$$

shall be called the "functionization" homomorphism.

One can describe this in a fancier way as follows. Let $\pi_\lambda : X_\infty \to X_\lambda$ denote the canonical projection induced from the projection $\prod_\lambda X_\lambda \to X_\lambda$. Consider the following commutative diagram (which follows from $\pi_\lambda = \pi_{\lambda\mu} \circ \pi_\mu (\lambda < \mu)$):

$$\begin{array}{c} F(X_\lambda) \\ \pi_{\lambda\mu}^* \downarrow \quad \searrow^{\pi_\lambda^*} \\ \quad\quad Fun(X_\infty, \mathbb{Z}) \\ \nearrow_{\pi_\mu^*} \\ F(X_\mu) \end{array}$$

Then it follows from a standard fact in the theory of inductive limits that the "functionization" homomorphism $\Psi : \varinjlim_{\lambda \in \Lambda} F(X_\lambda) \to Fun(X_\infty, \mathbb{Z})$ is nothing but the unique homomorphism such that the following diagram commutes:

$$\begin{array}{c} F(X_\lambda) \\ \rho^\lambda \swarrow \quad \searrow \pi_\lambda^* \\ F^{\mathrm{ind}}(X_\infty) \xrightarrow{\Psi} Fun(X_\infty, \mathbb{Z}). \end{array}$$

To avoid some possible confusion, the image $\Psi([\alpha_\lambda]) = \pi_\lambda^* \alpha_\lambda$ shall be denoted by $[\alpha_\lambda]_\infty$. For a constructible set $W_\lambda \in X_\lambda$, by the definition we have

$$[\mathbb{1}_{W_\lambda}]_\infty = \mathbb{1}_{\pi_\lambda^{-1}(W_\lambda)}.$$

$\pi_\lambda^{-1}(W_\lambda)$ is called a *proconstructible set (of level $\lambda$)* or a *cylinder set (of level $\lambda$)*, mimicking [Cr]. And the characteristic function supported on a proconstructible set (of level $\lambda$) is called a *procharacteristic function (of level $\lambda$)* and a finite linear combination of procharacteristic functions is called a *proconstructible function*. Let $F^{\mathrm{pro}}(X_\infty)$ denote the abelian group of all proconstructible functions on the proalgebraic variety $X_\infty = \varprojlim_{\lambda \in \Lambda} \{X_\lambda, \pi_{\lambda\mu} : X_\mu \to X_\lambda\}$. Thus we have the following

**Proposition 4.6.** *For a proalgebraic variety $X_\infty = \varprojlim_{\lambda \in \Lambda} \{X_\lambda, \pi_{\lambda\mu} : X_\mu \to X_\lambda\}$*

$$F^{\mathrm{pro}}(X_\infty) = \mathrm{Image}\left(\Psi : F^{\mathrm{ind}}(X_\infty) \to Fun(X_\infty, \mathbb{Z})\right) = \bigcup_\mu \pi_\mu^*(F(X_\mu)).$$

**Proposition 4.7.** *If the structure morphisms $\pi_{\lambda\mu} : X_\mu \to X_\lambda$ ($\lambda < \mu$) are all surjective, then for the proalgebraic variety $X_\infty = \varprojlim_{\lambda \in \Lambda} \{X_\lambda, \pi_{\lambda\mu} : X_\mu \to X_\lambda\}$ we have*

$$F^{\mathrm{ind}}(X_\infty) \cong F^{\mathrm{pro}}(X_\infty).$$

*Proof.* That all the structure morphisms $\pi_{\lambda\mu} : X_\mu \to X_\lambda$ ($\lambda < \mu$) are surjective implies that all the projections $\pi_\lambda : X_\infty \to X_\lambda$ are surjective. Which implies in turn that all the homomorphism $\pi_\lambda^* : F(X_\lambda) \to Fun(X_\infty, \mathbb{Z})$ are injective. Since the inductive limit is an exact functor, it follows that the "functionization" homomorphism $\Psi : \varinjlim_{\lambda \in \Lambda} F(X_\lambda) \to Fun(X_\infty, \mathbb{Z})$ is also injective. Thus we get the above isomorphism. $\square$

In the case of the arc space $\mathcal{L}(X)$ of a nonsingular variety $X$, since each structure morphism $\pi_{n(n+1)} : \mathcal{L}_{n+1}(X) \to \mathcal{L}_n(X)$ is always surjective, we get the following

**Corollary 4.8.** *For the arc space $\mathcal{L}(X)$ of a nonsingular variety $X$ we have the canonical isomorphism*

$$F^{\mathrm{ind}}(\mathcal{L}(X)) \cong F^{\mathrm{pro}}(\mathcal{L}(X)).$$

Suppose that $\Psi([\alpha_\mu]) = 0$, which means that $\Psi([\alpha_\mu])((x_\lambda) = \alpha_\mu(x_\mu) = 0$ for any $(x_\lambda) \in X_\infty$. Hence we have
$$\alpha_\mu\bigl(\pi_\mu(X_\infty)\bigr) = 0.$$
At the moment we do not know whether we can conclude $[\alpha_\mu] = 0$ from this condition. There is a very simple example such that $\alpha_\mu\bigl(\pi_\mu(X_\infty)\bigr) = 0, \pi_\mu(X_\infty) \neq X_\mu$ and $\alpha_\mu \neq 0$, but $[\alpha_\mu] = 0$ : Let $X_1 = \{a, b\}$ be a space of two different points, and let $X_n = \{a\}$ for any $n > 1$. Let $\pi_{12} : X_2 \to X_1$ be the injection map sending $a$ to $a$ and the other structure morphism $\pi_{n(n+1)} : X_{n+1} \to X_n$ is the identity for $n > 1$. Then the projective limit $X_\infty = \{(a)\}$ consists of one point $(a, a, a, \cdots)$. Let $\alpha_1 = p \cdot 1\!\!1_b \in F(X_1)$. Then we have $\alpha_1\left(\pi_1(X_\infty)\right) = 0, \pi_1(X_\infty) \neq X_1$ and $\alpha_1 \neq 0$, but $[\alpha_1] = 0$. We suspect that in general the "functionization" homomorphism $\Psi$ might be not necessarily injective, but we have not been able to find such an example yet:

**Question 4.9.** *Is the homomorphism* $\Psi : F^{\mathrm{ind}}(X_\infty) \to Fun(X_\infty, \mathbb{Z})$ *always injective ?*

Note that if we consider the topological situation and consider all functions, then the answer is certainly negative; e.g., consider the case of a decreasing sequence of subsets $X_n$ such that $\cap_{n=1}^\infty X_n = \emptyset$, in which case $F(X_\infty) = \{0\}$ by definition, but the inductive limit $\varinjlim_n F(X_n)$ contains all constant functions.

**Corollary 4.10.** *When $X$ is a nonsingular variety of dimension d, we have the following canonical Grothendieck class homomorphism*
$$\Gamma^{\mathrm{ind}} : F^{\mathrm{pro}}(\mathcal{L}(X)) \to K_0(\mathcal{V})_{[\mathbb{L}^d]}$$
*described by*
$$\Gamma^{\mathrm{ind}}\left([\alpha_n]_\infty\right) = \frac{\Gamma(\alpha_n)}{[\mathbb{L}]^{nd}}.$$
*In particular, we get that* $\Gamma^{\mathrm{ind}}\left(1\!\!1_{\mathcal{L}(X)}\right) = \Gamma^{\mathrm{ind}}([1\!\!1_X]_\infty) = [X]$.

Note that in the case of arc space $\mathcal{L}(X)$, since $\mathcal{L}_0(X) = X$, the indexed set is not $\mathbb{N}$ but $\{0\} \cup \mathbb{N}$. Hence the canonical one is not $\Gamma^{\mathrm{ind}}\left([\alpha_n]_\infty\right) = \dfrac{\Gamma(\alpha_n)}{[\mathbb{L}]^{(n-1)d}}$.

If $X$ is singular, the arc space $\mathcal{L}(X)$ is *not* the projective limit of a projective system of Zariski locally trivial fiber bundles with fiber being $\mathbb{C}^{\dim X}$ any longer and each projection morphism $\pi_{n(n+1)} : \mathcal{L}_{n+1}(X) \to \mathcal{L}_n(X)$ is complicated and thus as a proalgebraic variety $\mathcal{L}(X)$ is complicated. A crucial ingredient in studing motivic measure or motivic integration is the so-called *stable set* of the arc space $\mathcal{L}(X)$. A subset $A$ of the arc space $\mathcal{L}(X)$ is called *a stable set* if it is a cylinder set, i.e., $A = \pi_n^{-1}(C_n)$ for a constructible set $C_n$ in the $n$-th arc space $\mathcal{L}_n(X)$, such that the restriction of each projection $\pi_{m(m+1)}|_{\pi_{m+1}(A)} : \pi_{m+1}(A) \to \pi_m(A)$ for each $m \geq n$ is a Zariski locally fiber bundle with the fiber being $\mathbb{C}^{\dim X}$. So, our $\Gamma$-stable indconstructible function is a generalization of the characteristic function supported on this stable set.

Therefore we can see that our proalgebraic Grothendieck class homomorphism $\widetilde{\Gamma_G^{\mathrm{st.ind}}} : F_G^{\mathrm{st.ind}}(X_\infty) \to K_0(\mathcal{V})_G$ given in Theorem 4.3 (ii) is a generalization of the so-called motivic measure.

Before finishing this section, we give some remarks about another non-trivial and interesting generalized Euler–Poincaré characteristic, which is the *Hodge polynomial* (sometimes called the *Deligne–Hodge polynomial*, *E–polynomial* or *E–function*) defined via the theory of mixed Hodge structures [De1, De2] (e.g., see [Cr], [DK], [DL1, DL2], [Ve]). The existence of such a polynomial had been conjectured by J.-P. Serre before the theory of mixed Hodge structures was introduced by P. Deligne (see [F2, §4.5 and Notes to Chapter 4]).

For a complex algebraic variety $V$, we set

$$e^{p,q}(V) := \sum_{i \geq 0} (-1)^i h^{p,q}\left(H_c^i(V;\mathbb{C})\right),$$

which is called the *(p,q)-Hodge number* of the variety $X$ and the *Hodge polynomial* of $V$ is defined by

$$H_{u,v}(V) := \sum_{p,q} e^{p,q}(V) u^p v^q \in \mathbb{Z}[u,v].$$

The Hodge polynomial satisfies the following (e.g., see [DK]):

(i) If a complex variety $X$ has a finite stratification $X = \sqcup_i X_i$ by locally closed subvarieties $X_i$, then

$$H_{u,v}(X) = \sum_i H_{u,v}(X_i).$$

(ii) The Hodge polynomial is multiplicative; for complex algebraic varieties $X$ and $Y$

$$H_{u,v}(X \times Y) = H_{u,v}(X) \cdot H_{u,v}(Y).$$

(iii) (a more general version of (ii)) For a Zariski locally trivial fiber bundle $X \to Z$ with a fiber $F$

$$H_{u,v}(X) = H_{u,v}(F) \cdot H_{u,v}(Z).$$

A crucial difference between the topological Euler–Poincaré characteristic $\chi$ and the Hodge polynomial $H_{u,v}$ is that $H_{u,v}(X) = 0$ if and only if $X = \emptyset$. (Note that the degree of $H_{u,v}(X)$ is always $2dim_{\mathbb{C}} X$.) Another crucial difference is the property (iii); as to the topological Euler–Poincaré characteristic $\chi$ the equality $\chi(X) = \chi(F) \cdot \chi(Z)$ holds even for a *topological fiber bundle* $X \to Z$ with a fiber $F$.

It is obvious that one can define the following homomorphism

$$H_{u,v} : F(X) \to \mathbb{Z}[u,v]$$

defined by $H_{u,v}(\mathbb{1}_S) := H_{u,v}(S)$. With this definition, a naïve question is whether one can extend the Hodge polynomial homomorphism $H_{u,v} : F(X) \to \mathbb{Z}[u,v]$ to proalgebraic varieties as done above. If we think of the proof given in Example (3.6), one natural question is whether this Hodge polynomial is compatible with the pushforward $f_* : F(X) \to F(Y)$, and furthermore whether it can be extened to a higher homology class version, i.e., a natural transformation from the constructible function functor to a certain homology theory such that the homology group of a point equals $\mathbb{Z}[u,v]$. For an arbitrary one $(u,v)$, the answers for this question is negative becasuse of the above second difference (also, see [Jo, remarks after Proposition 3.17]). However, *in the special case when $(u,v) = (-y,1)$, for a nonsingular variety $X$ the Hodge polynomial $H_{-y,1}(X) \in \mathbb{Z}[y]$ is equal to the $\chi_y(X)$, the Hirzebruch $\chi_y$-genus of $X$* (see [Hi] and [HBJ]). Note that $H_{1,1}(X) = \chi(X)$ even if $X$ is singular and that its higher homology class version is the Chern–Schwartz–MacPherson class, and that for the special value $y = 0$ the Hirzebruch $\chi_0$-genus is the arithmetic genus and it has a higher homology class version for possibly singular varieties, i.e., Baum–Fulton–MacPherson's Riemann–Roch $\tau_* : \mathbf{K}_0 \to H_{*\mathbb{Q}}$ constructed in [BFM] and that for the special value $y = 1$ the Hirzebruch $\chi_1$-genus is the Thom–Hirzebruch $L$-class and it has a higher homology class version for possibly singular varieties, i.e., Cappell–Shaneson's homology $L$-class [CS1, CS2].

A (general) unified theory of characteristic classes of singular varieties or a theory unifying at least the above three characteristic classes of singular varieties has been looked for (e.g., see [Mac2] and [Y1]). And as above the reappearance of these three genera related to the corresponding characteristic classes through the theory of mixed Hodge structures was a kind of support to believe that there must be a reasonable positive solution to the above question in the special case when $(u,v) = (-y,1)$, in particular, with the speculation or simple-minded guess that Saito's mixed Hodge modules [Sai] would be a key to such a solution. And it turns out that one can give a positive solution to this question. For details of such a solution and many other related results, see [BSY2].

## 5. CHARACTERISTIC CLASSES OF PROALGEBRAIC VARIETIES

In this section we first consider generalizing MacPherson's Chern class transformation $c_* : F(X) \to H_*(X)$ to a category of proalgebraic varieties and modeled on this construction we consider general characteristic classes of proalgebraic varieties.

First we consider the infinite countable product $X_\infty := X^{\mathbb{N}}$ of a complex algebraic variety $X$ as a simple model case.

Let $X^n$ denote the Cartesian product of $n$ copies of the variety $X$. For each projection $\pi_{n(n+1)} : X^{n+1} \to X^n$ (projecting to the first $n$ factors), the pullback homomorphism $\pi_{n(n+1)}^* : F(X^n) \to F(X^{n+1})$ is the multiplication by the characteristic function $\mathbb{1}_X$ of the last factor $X$, i.e.,

$$\pi_{n(n+1)}^*(\alpha) = \alpha \times \mathbb{1}_X,$$

where $(\alpha \times \mathbb{1}_X)(y, x) := \alpha(y) \mathbb{1}_X(x) = \alpha(y)$. Then, using the cross product formula $c_*(\delta \times \omega) = c_*(\delta) \times c_*(\omega)$ of MacPherson's Chern class transformation $c_*$, due to M. Kwieciński [Kw] (cf. [KY]), we get the following commutative diagram

$$\begin{array}{ccc} F(X^n) & \xrightarrow{c_*} & H_*(X^n) \\ \times \mathbb{1}_X \downarrow & & \downarrow \times c_*(X) \\ F(X^{n+1}) & \xrightarrow{c_*} & H_*(X^{n+1}). \end{array}$$

So, if we set

$$H_{**}^{\mathrm{ind}}(X^{\mathbb{N}}) := \varinjlim_n \Big\{ \times c_*(X) : H_*(X^n) \to H_*(X^{n+1}) \Big\},$$

then we have a *proalgebraic MacPherson's Chern class homomorphism*:

$$c_*^{\mathrm{ind}} : F^{\mathrm{ind}}(X^{\mathbb{N}}) \to H_{**}^{\mathrm{ind}}(X^{\mathbb{N}}).$$

The proalgebraic Chern–Schwartz–MacPherson class $c_*^{\mathrm{ind}}(X_\infty) := c_*^{\mathrm{ind}}([\mathbb{1}_X]) = [c_*(X)]$.

Since we have

$$\begin{array}{ccc} H_*(X^n) & \xrightarrow{\int_{X_n}} & \mathbb{Z} \\ \times c_*(X) \downarrow & & \downarrow \times \chi(X) \\ H_*(X^{n+1}) & \xrightarrow{\int_{X_{n+1}}} & \mathbb{Z}, \end{array}$$

if we assume that $\chi(X) \neq 0$, then we get the proalgebraic integration

$$\int^{\mathrm{ind}} : H_{**}^{\mathrm{ind}}(X^{\mathbb{N}}) \to \mathbb{Z}\left[\frac{1}{\chi(X)}\right]$$

defined by $\displaystyle\int^{\mathrm{ind}} ([x_n]) = \left[\int_{X_n} x_n\right] = \frac{\int_{X_n} x_n}{\chi(X)^{n-1}}$, where $x_n \in H_*(X_n)$. And we also get that $\chi^{\mathrm{ind}} = \int^{\mathrm{ind}} \circ \, c_*^{\mathrm{ind}}$, which is a proalgebraic analogue of $\chi = \int \circ \, c_*$.

Second we consider the case of a proalgebraic variety $X_\infty = \varprojlim_{\lambda \in \Lambda} \Big\{ X_\lambda, \pi_{\lambda\mu} : X_\mu \to X_\lambda \Big\}$ of a projective system of smooth morphisms $\pi_{\lambda\mu} : X_\mu \to X_\lambda$. Here we recall the following *Verdier–Riemann–Roch formula for Chern class (abbr. VRR–Chern)* (see [FM], [Schü1] and [Y2]):

**Theorem 5.1.** *For a smooth morphism $f : X \to Y$ of possibly singular varieties $X$ and $Y$ the following diagram commutes:*

$$\begin{array}{ccc} F(Y) & \xrightarrow{c_*} & H_*(Y) \\ f^* \downarrow & & \downarrow c(T_f) \cap f^* \\ F(X) & \xrightarrow{c_*} & H_*(X). \end{array}$$

*Here $T_f$ is the relative tangent bundle of the smooth morphism $f$ and $f^* : H_*(Y) \to H_*(X)$ is the Gysin pullback homomorphism.*

From this **VRR–Chern** we get the following

**Corollary 5.2.** *(i) For a proalgebraic variety $X_\infty = \varprojlim_{\lambda \in \Lambda} \left\{ X_\lambda, \pi_{\lambda\mu} : X_\mu \to X_\lambda \right\}$ of a projective system of smooth morphisms $\pi_{\lambda\mu} : X_\mu \to X_\lambda$ we get the following proalgebraic MacPherson's Chern class homomorphism*

$$c_*^{\mathrm{ind}} : F^{\mathrm{ind}}(X_\infty) \to H_{**}^{\mathrm{ind}}(X_\infty) := \varinjlim_{\lambda \in \Lambda} \left\{ H_*(X_\lambda), c(T_{\pi_{\lambda\mu}}) \cap \pi_{\lambda\mu}{}^* : H_*(X_\lambda) \to H_*(X_\mu) \right\}.$$

*(ii) If the Euler–Poincaré characteristic of the fiber of each smooth morphism $\pi_{\lambda\mu}$ is non-zero and denoted by $e_{\lambda\mu}$, then we have the proalgebraic integration*

$$\int^{\mathrm{ind}} : H_{**}^{\mathrm{ind}}(X_\infty) \to \varinjlim_{\lambda \in \Lambda} \left\{ \times e_{\lambda\mu} : \mathbb{Z} \to \mathbb{Z} \right\}$$

*and we have the commutative diagram*

$$\begin{array}{ccc} F^{\mathrm{ind}}(X_\infty) & \xrightarrow{c_*^{\mathrm{ind}}} & H_{**}^{\mathrm{ind}}(X_\infty) \\ & \searrow \chi^{\mathrm{ind}} \quad \swarrow \int^{\mathrm{ind}} & \\ & \varinjlim_{\lambda \in \Lambda} \left\{ \times e_{\lambda\mu} : \mathbb{Z} \to \mathbb{Z} \right\}. & \end{array}$$

In this case we have that $c_*^{\mathrm{ind}}(X_\infty) = [c_*(X_\lambda)]$.

In order to generalize these results furthermore and also to capture the above $c_*^{\mathrm{ind}}$ as *a natural transformation*, we need to appeal to Fulton–MacPherson's bivariant homology theory [FM] and the Brasselet's bivariant Chern class [Br].

First, as we promised in §3, we discuss functorialities of $\mathcal{F}^{\mathrm{ind}}$ of a general bifunctor $\mathcal{F}$. Let $\{f_\lambda : X_\lambda \to Y_\lambda\}$ be a pro-morphism of pro-algebraic varieties $\{X_\lambda, \pi_{\lambda\mu} : X_\mu \to X_\lambda\}$ and $\left\{ Y_\lambda, \rho_{\lambda\mu} : Y_\mu \to Y_\lambda \right\}$. Then it follows from the contravariance of the bifunctor $\mathcal{F}$ that the following diagram commutes

$$\begin{array}{ccc} \mathcal{F}(Y_\lambda) & \xrightarrow{f_\lambda^*} & \mathcal{F}(X_\lambda) \\ \rho_{\lambda\mu}{}^* \downarrow & & \downarrow \pi_{\lambda\mu}{}^* \\ \mathcal{F}(Y_\mu) & \xrightarrow{f_\mu^*} & \mathcal{F}(X_\mu), \end{array}$$

which in turn implies that the pullback homomorphism $f_\infty^* := \varinjlim\{f_\lambda^*\} : \mathcal{F}^{\mathrm{ind}}(Y_\infty) \to \mathcal{F}^{\mathrm{ind}}(X_\infty)$ is a contravariantly functorial. However, to claim the covariance of $\mathcal{F}^{\mathrm{ind}}$, we need the following requirements; one for the bifunctor $\mathcal{F}$ and one for the pro-morphism $\{f_\lambda : Y_\lambda \to X_\lambda\}_{\lambda \in \Lambda}$:

**Definition 5.3.** If a bifunctor $\mathcal{F} : \mathcal{V} \to \mathcal{A}$ satisfies the following two properties (M-1) and (M-2), then it is called a *Mackey functor*:

(M-1): for any fiber square in $\mathcal{V}$

$$\begin{array}{ccc} X' & \xrightarrow{g'} & X \\ f' \downarrow & & \downarrow f \\ Y' & \xrightarrow{g} & Y \end{array}$$

the following diagram commutes

$$\begin{array}{ccc} \mathcal{F}(X) & \xrightarrow{g'^*} & \mathcal{F}(X') \\ f_* \downarrow & & \downarrow f'_* \\ \mathcal{F}(Y) & \xrightarrow{g^*} & \mathcal{F}(Y'), \end{array}$$

(M-2): $\qquad \mathcal{F}(X \coprod Y) = \mathcal{F}(X) \oplus \mathcal{F}(Y).$

The constructible function functor $F(X)$ and the relative Grothendieck group functor $K_0(\mathcal{V}/X)$ are Mackey functors.

The notion of Mackey functor was introduced by A. W. Dress [Dr1, Dr2] (also see [Bou] and [TW]) in the representation theory of finite groups. In what follows, the property we need is just the property (M-1), which is sometimes called the *base change formula* and a bifunctor satisfying (M-1) is called a *pre-Mackey functor*.

Let $\mathcal{F}, \mathcal{G} : \mathcal{V} \to \mathcal{A}$ be two (pre-)Mackey functors, and let $\Theta : \mathcal{F} \to \mathcal{G}$ be a natural transformation, i.e., for any morphism $f : X \to Y$ the following diagrams commute:

$$\begin{array}{ccc} \mathcal{F}(X) & \xrightarrow{\Theta_X} & \mathcal{G}(X) \\ \mathcal{F}_*(f) \downarrow & & \downarrow \mathcal{G}_*(f) \\ \mathcal{F}(Y) & \xrightarrow{\Theta_Y} & \mathcal{G}(Y) \end{array} \qquad \begin{array}{ccc} \mathcal{F}(Y) & \xrightarrow{\Theta_Y} & \mathcal{G}(Y) \\ \mathcal{F}^*(f) \downarrow & & \downarrow \mathcal{G}^*(f) \\ \mathcal{F}(X) & \xrightarrow{\Theta_X} & \mathcal{G}(X). \end{array}$$

From now on, unless some confusion is possible, we just denote $f_*$ for both $\mathcal{F}_*(f)$ and $\mathcal{G}_*(f)$, $f^*$ for both $\mathcal{F}^*(f)$ and $\mathcal{G}^*(f)$, and $\Theta$ for $\Theta_X, \Theta_Y$ without subscripts.

**Definition 5.4.** Definition (5.4)Let $\{f_\lambda : X_\lambda \to Y_\lambda\}_{\lambda \in \Lambda}$ be a pro-morphism of pro-algebraic varieties $\left\{X_\lambda, \pi_{\lambda\mu} : X_\mu \to X_\lambda\right\}$ and $\left\{Y_\lambda, \rho_{\lambda\mu} : Y_\mu \to Y_\lambda\right\}$. If the following commutative diagram for $\lambda < \mu$

$$\begin{array}{ccc} X_\mu & \xrightarrow{f_\mu} & Y_\mu \\ \pi_{\lambda\mu} \downarrow & & \downarrow \rho_{\lambda\mu} \\ X_\lambda & \xrightarrow{f_\lambda} & Y_\lambda \end{array}$$

is a fiber square, then we call the pro-morphism $\{f_\lambda : X_\lambda \to Y_\lambda\}_{\lambda \in \Lambda}$ a *fiber-square pro-morphism*, abusing words.

**Theorem 5.5.** *(i) Let $\mathcal{F} : \mathcal{V} \to \mathcal{A}$ be a (pre-)Mackey functor. Then for a fiber-square pro-morphism $\{f_\lambda : X_\lambda \to Y_\lambda\}_{\lambda \in \Lambda}$ of pro-algebraic varieties $\left\{X_\lambda, \pi_{\lambda\mu} : X_\mu \to X_\lambda\right\}$ and $\left\{Y_\lambda, \rho_{\lambda\mu} : Y_\mu \to Y_\lambda\right\}$, the pushforward homomorphism*

$$f_{\infty *} := \varinjlim_{\lambda \in \Lambda} \{f_{\lambda *}\} : \mathcal{F}^{\mathrm{ind}}(X_\infty) \to \mathcal{F}^{\mathrm{ind}}(Y_\infty)$$

*is covariantly functorial.*

(ii) *Furthermore, for a projective system* $P = \{p_{\lambda\mu}\}$ *of non-zero elements of* $\mathcal{R}$ *the following diagram commutes:*

$$\begin{array}{ccc}
\mathcal{F}_P^{\text{st.ind}}(X_\infty) & \xrightarrow{\chi_{\mathcal{F}}^{\text{ind}}} & \\
f_{\infty *} \downarrow & & \varinjlim \{\times p_{\lambda\mu} : \mathcal{R} \to \mathcal{R}\}. \\
\mathcal{F}_P^{\text{st.ind}}(Y_\infty) & \xrightarrow{\chi_{\mathcal{F}}^{\text{ind}}} & 
\end{array}$$

(iii) *Let* $\mathcal{F}, \mathcal{G} : \mathcal{V} \to \mathcal{A}$ *be two (pre-)Mackey functors and let* $\Theta : \mathcal{F} \to \mathcal{G}$ *be a natural transformation. For a projective system* $P = \{p_{\lambda\mu}\}$ *of non-zero elements* $p_{\lambda\mu}$ *of* $\mathcal{R}$ *and a fiber-square pro-morphism* $\{f_\lambda : X_\lambda \to Y_\lambda\}_{\lambda \in \Lambda}$ *of pro-algebraic varieties* $\{X_\lambda, \pi_{\lambda\mu} : X_\mu \to X_\lambda\}$ *and* $\{Y_\lambda, \rho_{\lambda\mu} : Y_\mu \to Y_\lambda\}$, *we have the commutative diagram*

$$\begin{array}{ccc}
\mathcal{F}^{\text{ind}}(X_\infty) & \xrightarrow{\Theta^{\text{ind}}} & \mathcal{G}^{\text{ind}}(X_\infty) \\
f_{\infty *} \downarrow & & \downarrow f_{\infty *} \\
\mathcal{F}^{\text{ind}}(Y_\infty) & \xrightarrow{\Theta^{\text{ind}}} & \mathcal{G}^{\text{ind}}(Y_\infty).
\end{array}$$

(iv) *Furthermore we suppose that* $\mathcal{F}(pt) = \mathcal{G}(pt) = \mathcal{R}$ *is a commutative ring with a unit and* $\Theta : \mathcal{R} = \mathcal{F}(pt) \to \mathcal{R} = \mathcal{G}(pt)$ *is the identity. Then we have the following commutative diagrams:*

$$\begin{array}{ccc}
\mathcal{F}_P^{\text{st.ind}}(X_\infty) & \xrightarrow{\Theta^{\text{ind}}} & \mathcal{G}_P^{\text{st.ind}}(X_\infty) \\
& \chi_{\mathcal{F}}^{\text{ind}} \searrow \quad \swarrow \chi_{\mathcal{G}}^{\text{ind}} & \\
& \varinjlim_{\lambda \in \Lambda} \{\times p_{\lambda\mu} : \mathcal{R} \to \mathcal{R}\}. & 
\end{array}$$

*Proof.* It suffices to see that for a fiber-square pro-morphism $\{f_\lambda : X_\lambda \to Y_\lambda\}_{\lambda \in \Lambda}$ of pro-algebraic varieties $\{X_\lambda, \pi_{\lambda\mu} : X_\mu \to X_\lambda\}$ and $\{Y_\lambda, \rho_{\lambda\mu} : Y_\mu \to Y_\lambda\}$, we get the following commutative cubic diagram:

$$\begin{array}{ccccc}
& & \mathcal{F}(X_\lambda) & \xrightarrow{f_{\lambda *}} & \mathcal{F}(Y_\lambda) \\
& \Theta \swarrow & & \Theta \swarrow & \\
\mathcal{G}(X_\lambda) & \xrightarrow{f_{\lambda *}} & \mathcal{G}(X_\lambda) & & \rho_{\lambda\mu}^* \downarrow \\
& & \pi_{\lambda\mu}^* \downarrow & \rho_{\lambda\mu}^* & \\
\pi_{\lambda\mu}^* \downarrow & & \mathcal{F}(X_\mu) & \xrightarrow{f_{\mu *}} & \mathcal{F}(Y_\mu) \\
& \Theta \swarrow & & \Theta \swarrow & \\
\mathcal{G}(X_\mu) & \xrightarrow{f_{\mu *}} & \mathcal{G}(Y_\mu). & & 
\end{array}$$

$\square$

The homology theory is not a (pre-)Mackey functor, and to get a generalized version of Corollary 5.2, we use Fulton–MacPherson's bivariant homology theory $\mathbb{H}$, constructed from the cohomology theory. For a morphism $f : X \to Y$, choose a morphism $\phi : X \to \mathbb{R}^n$ such that $\Phi := (f, \phi) : X \to Y \times \mathbb{R}^n$ is a closed embedding. Then the $i$-th bivariant homology group $\mathbb{H}^i(X \xrightarrow{f} Y)$ is defined by

$$\mathbb{H}^i(X \xrightarrow{f} Y) := H^{i+n}(Y \times \mathbb{R}^n, Y \times \mathbb{R}^n \setminus X_\phi),$$

where $X_\phi$ is defined to be the image of the morphism $\Phi = (f, \phi)$. The definition is independent of the choice of $\phi$. Note that instead of taking the Euclidean space $\mathbb{R}^n$ we can take a manifold $M$ so that $i : X \to M$ is a closed embedding and then consider the graph embedding $f \times i : X \to Y \times M$. See [FM, §3.1] for more details of $\mathbb{H}$. In particular, note that if $Y$ is a point $pt$, $\mathbb{H}(X \to pt)$ is isomorphic to the homology group $H_*(X)$ of the source variety $X$.

W. Fulton and R. MacPherson conjectured or posed as a question the existence of a so-called bivariant Chern class and J.-P. Brasselet [Br] solved it:

**Theorem 5.6.** *On the category of embeddable complex analytic varieties and cellular morphisms, there exists a Grothendieck transformation*

$$\gamma^{\mathrm{Br}} : \mathbb{F} \to \mathbb{H}$$

*satisfying the normalization condition that* $\gamma^{\mathrm{Br}}(\mathbb{1}_\pi) = c(TX) \cap [X]$ *for $X$ smooth, where* $\pi : X \to pt$ *and* $\mathbb{1}_\pi = \mathbb{1}_X$.

Note that for a morphism $\pi : X \to pt$ from a variety $X$ to a point $pt$, $\gamma^{\mathrm{Br}} : \mathbb{F}(X \to pt) \to \mathbb{H}(X \to pt)$ is nothing but the original MacPherson's Chern class transformation $c_* : F(X) \to H_*(X)$.

As observed in §3, we get the following

**Corollary 5.7.** *(A generalized Verdier–Riemann–Roch formula for Chern class) For a bivariant constructible function* $\alpha \in \mathbb{F}(X \xrightarrow{f} Y)$ *we have the following commutative diagram:*

$$\begin{array}{ccc}
F(Y) & \xrightarrow{c_*} & H_*(Y) \\
{\scriptstyle \alpha \bullet_\mathbb{F} = \alpha \cdot f^*} \downarrow & & \downarrow {\scriptstyle \gamma^{\mathrm{Br}}(\alpha) \bullet_\mathbb{H}} \\
F(X) & \xrightarrow{c_*} & H_*(X).
\end{array}$$

*In particular, for an Euler morphism we have the following commutative diagram:*

$$\begin{array}{ccc}
F(Y) & \xrightarrow{c_*} & H_*(Y) \\
{\scriptstyle \mathbb{1}_f \bullet_\mathbb{F} = f^*} \downarrow & & \downarrow {\scriptstyle \gamma^{\mathrm{Br}}(\mathbb{1}_f) \bullet_\mathbb{H}} \\
F(X) & \xrightarrow{c_*} & H_*(X).
\end{array}$$

For a more generalized Verdier–Riemann–Roch theorem for Chern class, see [Schü1]. The homomorphism $\gamma^{\mathrm{Br}}(\mathbb{1}_f) \bullet_\mathbb{H}$ shall be denoted by $f^{**}$.

Using Corollary 5.7, we get the following

**Theorem 5.8.** *(i) Let $\{f_\lambda : X_\lambda \to Y_\lambda\}$ be a fiber-square pro-morphism between two proalgebraic varieties $\left\{X_\lambda, \pi_{\lambda\mu} : X_\mu \to X_\lambda\right\}$ and $\left\{Y_\lambda, \rho_{\lambda\mu} : Y_\mu \to Y_\lambda\right\}$ with structure morphisms being Euler morphisms. Then we have the following commutative diagram:*

$$\begin{array}{ccc}
F^{\mathrm{ind}}(X_\infty) & \xrightarrow{c_*^{\mathrm{ind}}} & H_*^{\mathrm{ind}}(X_\infty; \{\gamma^{\mathrm{Br}}(\mathbb{1}_{\pi_{\lambda\mu}})\}) \\
{\scriptstyle f_{\infty *}} \downarrow & & \downarrow {\scriptstyle f_{\infty *}} \\
F^{\mathrm{ind}}(Y_\infty) & \xrightarrow{c_*^{\mathrm{ind}}} & H_*^{\mathrm{ind}}(Y_\infty; \{\gamma^{\mathrm{Br}}(\mathbb{1}_{\rho_{\lambda\mu}})\}).
\end{array}$$

*(ii) Let $X_\infty = \varprojlim_{\lambda \in \Lambda}\left\{X_\lambda, \pi_{\lambda\mu} : X_\mu \to X_\lambda\right\}$ be a proalgebraic variety such that for each $\lambda < \mu$ the structure morphism $\pi_{\lambda\mu} : X_\mu \to X_\lambda$ is an Euler proper morphism*

*(hence surjective) of topologically connected algebraic varieties with the constant Euler–Poincaré characteristic $\chi_{\lambda\mu}$ of the fiber of the morphism $\pi_{\lambda\mu}$ being non-zero. Then we get the commutative diagram*

$$\begin{array}{ccc} F^{\mathrm{ind}}(X_\infty) & \xrightarrow{c_*^{\mathrm{ind}}} & H_*^{\mathrm{ind}}(X_\infty; \{\gamma^{\mathrm{Br}}(\mathbb{1}_{\pi_{\lambda\mu}})\}) \\ & \searrow^{\chi^{\mathrm{ind}}} \quad \swarrow_{\int^{\mathrm{ind}}} & \\ & \varinjlim_\Lambda \{\times \chi_{\lambda\mu} : \mathbb{Z} \to \mathbb{Z}\} & \end{array}$$

*Proof.* As shown in the proof of Theorem 5.5, it suffices to show the commutativity of the following cubic diagram:

$$\begin{array}{ccc} F(X_\lambda) & \xrightarrow{f_{\lambda*}} & F(Y_\lambda) \\ \downarrow & & \downarrow \\ H_*(X_\lambda) & \xrightarrow{f_{\lambda*}} & H_*(Y_\lambda) \\ \downarrow \pi_{\lambda\mu}^{**} & & \downarrow \rho_{\lambda\mu}^{**} \\ F(X_\mu) & \xrightarrow{f_{\mu*}} & F(Y_\mu) \\ \downarrow & & \downarrow \\ H_*(X_\mu) & \xrightarrow{f_{\mu*}} & H_*(Y_\mu) \end{array}$$

The commutativity of the top and bottom squares is due to the naturality of MacPherson's Chern class transformation, the commutativity of the right and left squares is due to the above Corollary 5.7, and the commutativity of the square in the back is due to the Mackey property of the constructible function functor $F$. Thus it remains to see only the commutativity of the square in the front, i.e., for any $x \in H_*(X_\lambda)$

$$\rho_{\lambda\mu}^{**}\left(f_{\lambda*}(x)\right) = f_{\mu*}\left(\pi_{\lambda\mu}^{**}(x)\right),$$

which is more precisely

$$\gamma^{\mathrm{Br}}(\mathbb{1}_{\rho_{\lambda\mu}}) \bullet_{\mathbb{H}} f_{\lambda\star}(x) = f_{\mu*}\left(\gamma^{\mathrm{Br}}(\mathbb{1}_{\pi_{\lambda\mu}}) \bullet_{\mathbb{H}} x\right).$$

Since $\mathbb{1}_{\pi_{\lambda\mu}} = f_\lambda^\star \mathbb{1}_{\rho_{\lambda\mu}}$ and the Grothendieck transformation $\gamma^{\mathrm{Br}} : \mathbb{F} \to \mathbb{H}$ is compatible with the pullback operation, the above equality becomes

$$\gamma^{\mathrm{Br}}(\mathbb{1}_{\rho_{\lambda\mu}}) \bullet_{\mathbb{H}} f_{\lambda\star}(x) = f_{\mu*}\left(f_\lambda^\star \gamma^{\mathrm{Br}}(\mathbb{1}_{\rho_{\lambda\mu}}) \bullet_{\mathbb{H}} x\right).$$

And it turns out that this equality is nothing but *the projection formula* of the Bivariant Theory [FM, §2.2, $(A_{123})$] for the following diagram and for the bivariant homology theory $\mathbb{H}$:

$$\begin{array}{ccc} X_\mu & \xrightarrow{f_\mu} & Y_\mu \\ \pi_{\lambda\mu} \downarrow & & \downarrow \rho_{\lambda\mu} \\ X_\lambda & \xrightarrow{f_\lambda} & Y_\lambda \longrightarrow pt. \end{array}$$

Thus we get the theorem. $\square$

Following the above construction, similarly we can get a proalgebraic version of Baum–Fulton–MacPherson's Riemann–Roch $\tau_* : \mathbf{K}_0 \to H_{*\mathbb{Q}}$ constructed in [BFM], using the bivariant Riemann–Roch theorem ([Fu] and [FM]). And a much more general theorem is the following *characteristic classes of proalgebraic varieties*:

**Theorem 5.9.** *(i) Let $\gamma : \mathbb{B} \to \mathbb{B}'$ be a Grothendieck transformation between two bivariant theories $\mathbb{B}, \mathbb{B}' : \mathcal{C} \to \mathcal{A}$ and let $\{(\pi_{\lambda\mu}; b_{\lambda\mu}) : X_\mu \to X_\lambda\}$ be a projective system of bivariant-class-equipped morphisms. Then we get the following pro-version of the natural transformation $\gamma_* : \mathbb{B}_* \to \mathbb{B}'_*$:*

$$\gamma_*^{\mathrm{ind}} : \mathbb{B}_*^{\mathrm{ind}}\Big(X_\infty; \{b_{\lambda\mu}\}\Big) \to \mathbb{B}'_*{}^{\mathrm{ind}}\Big(X_\infty; \{\gamma(b_{\lambda\mu})\}\Big).$$

*(ii) Let $\{f_\lambda : X_\lambda \to Y_\lambda\}$ be a fiber-square pro-morphism between two projective systems $\{(\pi_{\lambda\mu}; b_{\lambda\mu}) : X_\mu \to X_\lambda\}$ and $\{(\rho_{\lambda\mu}; d_{\lambda\mu}) : Y_\mu \to Y_\lambda\}$ of bivariant-class-equipped morphisms such that $b_{\lambda\mu} = f_\lambda^\star d_{\lambda\mu}$. Then we have the following commutative diagram:*

$$\begin{array}{ccc}
\mathbb{B}_*^{\mathrm{ind}}(X_\infty; \{b_{\lambda\mu}\}) & \xrightarrow{\gamma_*^{\mathrm{ind}}} & \mathbb{B}'_*{}^{\mathrm{ind}}(X_\infty; \{\gamma(b_{\lambda\mu})\}) \\
f_{\infty *} \downarrow & & \downarrow f_{\infty *} \\
\mathbb{B}_*^{\mathrm{ind}}(Y_\infty; \{d_{\lambda\mu}\}) & \xrightarrow{\gamma_*^{\mathrm{ind}}} & \mathbb{B}'_*{}^{\mathrm{ind}}(Y_\infty; \{\gamma(d_{\lambda\mu})\}).
\end{array}$$

*(iii) Let $\mathbb{B}_*(pt) = \mathbb{B}'_*(pt)$ be a commutative ring $\mathcal{R}$ with a unit and we assume that the homomorphism $\gamma : \mathbb{B}_*(pt) \to \mathbb{B}'_*(pt)$ is the identity. Let $P = \{p_{\lambda\mu}\}$ be a projective system of non-zero elements $p_{\lambda\mu} \in \mathcal{R}$. Then we get the commutative diagram*

$$\mathbb{B}_{*,P}^{\mathrm{st.ind}}\Big(X_\infty; \{b_{\lambda\mu}\}\Big) \xrightarrow{\gamma_*^{\mathrm{ind}}} \mathbb{B}'_{*,P}{}^{\mathrm{st.ind}}\Big(X_\infty; \{\gamma(b_{\lambda\mu})\}\Big)$$

with maps $\chi_{\mathbb{B}_*}^{\mathrm{ind}}$ and $\chi_{\mathbb{B}'_*}^{\mathrm{ind}}$ to $\varinjlim_{\lambda \in \Lambda}\Big\{\times p_{\lambda\mu} : \mathcal{R} \to \mathcal{R}\Big\}.$

*Proof.* As in Theorem 5.8, it follows from the following commutative diagram:

$$\begin{array}{c}\text{(cube diagram involving } \mathbb{B}_*(X_\lambda), \mathbb{B}_*(Y_\lambda), \mathbb{B}'_*(X_\lambda), \mathbb{B}'_*(Y_\lambda), \mathbb{B}_*(X_\mu), \mathbb{B}_*(Y_\mu), \mathbb{B}'_*(X_\mu), \mathbb{B}'_*(Y_\mu) \\ \text{with maps } f_{\lambda*}, f_{\mu*}, \gamma_*, b_{\lambda\mu}\bullet, d_{\lambda\mu}\bullet, \gamma(b_{\lambda\mu})\bullet, \gamma(d_{\lambda\mu})\bullet)\end{array}$$

□

**Remark 5.10.** (i) As shown in [BSY1] (also see [EY1, EY2], [Schü2], [Y3, Y4, Y5, Y6]), a natural transformation bewteen two covariant functors commuting with exterior products is always extended to a Grothendieck transformation between their associated bivariant theories. Hence, as done in this section, it follows that such a natural transformation between two covaraint functors can be extended to a natural transformations between the proalgebraic versions of the covaraint functors for the category of proalgebraic varieties.
(ii) A much more abstract situation dealing with bifunctors is treated in [Y8].

## 6. GREEN FUNCTORS AND GROTHENDIECK–GREEN FUNCTORS

In this section we discuss a uniqueness of the canonical homomorphism $e : K_0(\mathcal{V}/X) \to F(X)$ defined by $e([Y \xrightarrow{f} X]) := f_* \mathbb{1}_Y$. A good reference for this section is [Bou].

**Definition 6.1.** (Green functors) A Green functor $G = (G^*, G_*)$ is a Mackey functor endowed with a bilinear map (or an exterior product)
$$G(X) \times G(Y) \to G(X \times Y)$$
denoted by $(x, y) \mapsto x \times y$ which are bifunctorial, associative and unitary, in the following sense:

(G-I) (bifunctoriality) for morphisms $f : X \to X'$ and $g : Y \to Y'$ the following diagrams commute:

$$\begin{array}{ccc} G(X) \times G(Y) & \xrightarrow{\times} & G(X \times Y) \\ f_* \times g_* \downarrow & & \downarrow (f \times g)_* \\ G(X') \times G(Y') & \xrightarrow{\times} & G(X' \times Y'), \end{array} \qquad \begin{array}{ccc} G(X') \times G(Y') & \xrightarrow{\times} & G(X' \times Y') \\ f^* \times g^* \downarrow & & \downarrow (f \times g)^* \\ G(X) \times G(Y) & \xrightarrow{\times} & G(X \times Y). \end{array}$$

(G-II) (associativity) $(x \times y) \times z = x \times (y \times z)$ for $x \in G(X), y \in G(Y), z \in G(Z)$. To be more precise, the following square

$$\begin{array}{ccc} G(X) \times G(Y) \times G(Z) & \xrightarrow{Id_{G(X)} \times (\times)} & G(X) \times G(Y \times Z) \\ (\times) \times Id_{G(Z)} \downarrow & & \downarrow \times \\ G(X \times Y) \times G(Z) & \xrightarrow{\times} & G(X \times Y \times Z) \end{array}$$

is commutative, up to identifications $(X \times Y) \times Z \cong X \times Y \times Z \cong X \times (Y \times Z)$.

(G-III) (unitarity) For a point $pt$ there exists a unit $1_G \in G(pt)$ such that for any $x \in G(X)$
$$p_{1*}(x \times 1_G) = x = p_{2*}(1_G \times x).$$
Here $p_1 : X \times pt \to X$ and $p_2 : pt \times X \to X$ are the projections (which are in fact isomorphisms).

The corresponding ones in the representations of finite groups is called the Burnside ring or the Burnside functor (e.g., see [Bou]).

**Remark 6.2.** For a Green functor $G$, by the identification $pt \times pt \cong pt$, the abelian group $G(pt)$ becomes a ring with the exterior product operation and the other abelian group $G(X)$ is a $G(pt)$-module.

**Remark 6.3.** The theory of *algebraic cobordism* introduced by M. Levine and F. Morel [LM] is a much finer theory of Green functors in the following sense. The pushforward homomorphisms are considered only for projective morphisms and the pullback homomorphisms are considered only for smooth morphisms. In such a restricted situation, it shall be called *a restricted Green functor*. Such a theory is sometimes called *a Borel–Moore functor with products* (e.g., see [LP]). Furthermore, if it is required that it is "oriented", i.e., it is equipped with Chern class operations, then such a theory is called *an oriented Borel–Moore functor with products* (e.g., see [LP]), and so it may also be called *an oriented restricted Green functor*.

The constructible function functor $F(X)$ is a Green functor if we consider the exterior product
$$\times : F(X) \times F(Y) \to F(X \times Y), \quad (\alpha, \beta) \mapsto \alpha \times \beta$$
defined by
$$(\alpha \times \beta)(x, y) := \alpha(x)\beta(y).$$
And the relative Grothendieck group $K_0(\mathcal{V}/X)$ is also a Green functor, if we consider the exterior product
$$K_0(\mathcal{V}/X) \times K_0(\mathcal{V}/Y) \to K_0(\mathcal{V}/X \times Y)$$

defined by the product of morphisms:

$$[X' \xrightarrow{f} X] \times [Y' \xrightarrow{g} Y] := [X' \times Y' \xrightarrow{f \times g} Y].$$

If $G, G'$ are Green functors on a category $\mathcal{C}$, a morphism or a natural transformation $\tau$ from $G$ to $G'$ is a natural transformation of Mackey functors $G$ and $G'$ which is compatible with exterior products, i.e., such that for a variety $X$ the following diagram commutes:

$$\begin{array}{ccc} G(X) \times G(Y) & \xrightarrow{\times} & G(X \times Y) \\ \tau_X \times \tau_Y \downarrow & & \downarrow \tau_{X \times Y} \\ G'(X) \times G'(Y) & \xrightarrow{\times} & G'(X \times Y). \end{array}$$

If moreover $\tau_{pt} : G(pt) \to G'(pt)$ sends the unit to the unit, then the natural transformation $\tau$ is called *unitary*.

**Definition 6.4.** If a Green functor $G = (G^*, G_*)$ satisfies the following "additivity": for a closed subvariety $Z \subset Y$

$$p_Y^*(1_G) = i_{Y-Z*}i_{Y-Z}^*p_Y^*(1_G) + i_{Z*}i_Z^*p_Y^*(1_G),$$

then it is called a *Grothendieck–Green* functor. Here we let $p_W : X \to pt$ be the map to a point for a variety $W$.

The constructible function functor $F(X)$ and the relative Grothendieck group functor $K_0(\mathcal{V}/X)$ are both Grothendieck–Green functors. Another highly nontrivial example of a Grothendieck–Green functor is the Grothendieck ring $K_0\Big(D^b(MHM(X))\Big)$ of the derived category of mixed Hodge modules with the natural t-structure (see [Getz, Proposition 3.9 and Definition 4.3]).

**Remark 6.5.** In the definition of the relative Grothendieck group $K_0(\mathcal{V}/X)$, an element is the isomorphism class $[Y \xrightarrow{f} X]$ of a morphism $f : Y \to X$. So, one might be tempted to include the following requirement in the above definition of the Grothendieck–Green functor: $h_*h^*p_Y^*(1_G) = p_Y^*(1_G)$ for an isomorphism $h : Y' \cong Y$. But it turns out that this requirement automatically follows from the Mackey property: indeed, consider the following fiber square

$$\begin{array}{ccc} Y' & \xrightarrow{h} & Y \\ h \downarrow & & \downarrow id_Y \\ Y & \xrightarrow{id_Y} & Y, \end{array}$$

from which we get that $h_*h^* = id_{G(Y)}$.

The following theorem is an algebro-geometric analogue of [Bou, Proposition 2.4.4]):

**Theorem 6.6.** *(A fundamental principle for Grothendieck–Green functors) For any unitary Grothendieck–Green functor $G : \mathcal{V} \to \mathcal{A}$, there exists a unique unitary natural transformation of Grothendieck–Green functors*

$$\tau : K_0(\mathcal{V}/\ ) \to G.$$

*Proof.* Let $[W \xrightarrow{h} X] \in K_0(\mathcal{V}/X)$ and let $p_W : W \to pt$ be the map to a point and let $i_W : W \to X$ be the inclusion. Then $[W \xrightarrow{h} X]$ can be expressed as

$$[W \xrightarrow{h} X] = h_*p_W^*([pt \to pt]).$$

Let $G$ be another Grothendieck–Green functor. If there exists a unitary natural transformation $\tau : K_0(\mathcal{V}/\ ) \to G$, then it follows from the naturality and unitarity that we have to have

$$\tau_X([W \xrightarrow{h} X]) = \tau_X\left(h_* p_W^*([pt \to pt])\right)$$
$$= h_* p_W^*(1_G)$$

So, all we have to do is to show that

$$\tau_X([W \xrightarrow{h} X]) := h_* p_W^*(1_G)$$

gives us a natural transformation between two Grothendieck–Green functors, and then we are done. Since the proof is straightforward, it is left for the reader. □

As a corollary of this theorem, a unitary natural transformation from $e : K_0(\mathcal{V}/X) \to F(X)$ *has to* be defined by $e([Y \xrightarrow{h} X]) := f_* \mathbb{1}_Y$.

**Remark 6.7.** In the above theorem, one cannot replace the Grothendieck–Green functor $K_0(\mathcal{V}/\ )$ by the constructible function Grothendieck–Green functor $F$. For the characteristic function $\mathbb{1}_W \in F(X)$ for a subvariety $W \subset X$ we have that, as in the above proof, $\mathbb{1}_W$ can be expressed as $\mathbb{1}_W = i_{W*} p_W^*(\mathbb{1}_{pt})$, where $i_W : W \to X$ be the inclusion. Hence, as in the above proof, we could define $\tau_X(\mathbb{1}_W) := (i_W)_* p_W^*(1_G)$. Then, all the arguments of the above proof perfectly work even for the constructible function Grothendieck–Green functor $F$, except for the naturality of the pushforward:

$$\begin{array}{ccc} F(X) & \xrightarrow{\tau_X} & G(X) \\ f_* \downarrow & & \downarrow G(f)_* \\ F(Y) & \xrightarrow{\tau_Y} & G(Y). \end{array}$$

In fact, one can see that this does not already hold for $G = K_0(\mathcal{V}/\ )$. Indeed, if it were the case, the uniqueness of such a unitary natural transformation would imply that for any variety $X$ we should have the isomorphism $K_0(\mathcal{V}/X) \cong F(X)$ and hence, in particular, we would have the isomorphism $K_0(\mathcal{V}/pt) \cong F(pt) \cong \mathbb{Z}$, which contradicts the recent result of Poonen [Po] that the Grothendieck ring $K_0(\mathcal{V})$ of varieties over a field of characteristic zero is not a domain.

**Remark 6.8.** Surely $i_F : F(X) \to K_0(\mathcal{V}/X)$ defined by $i_F(\mathbb{1}_W) := [W \xrightarrow{i_W} X]$ is injective. However, the above remark implies that this injective transformation cannot be a unitary natural transformation between the two Grothendieck–Green functors.

**Remark 6.9.** Applying the above natural transformation $e : K_0(\mathcal{V}/?) \to F(?)$ to a morphism $X \to pt$, we get the commutative diagram

$$\begin{array}{ccc} K_0(\mathcal{V}/X) & \xrightarrow{\chi_{\mathrm{Gro}}} & K_0(\mathcal{V}) \\ e \downarrow & & \downarrow e \\ F(X) & \xrightarrow{\chi} & \mathbb{Z}. \end{array}$$

And as we observed in §2, MacPherson's Chern class transformation $c_* : F(X) \to H_*(X)$ is a higher homology class version of $\chi : F(X) \to \mathbb{Z}$. It turns out that we can get a similar result for the above homomorphism $\chi_{\mathrm{Gro}} : K_0(\mathcal{V}/X) \to K_0(\mathcal{V})$ in such a way that it fits in the above commutative diagram; namely we can show the existence of a reasonable abelian group $Ab(\mathcal{V}/X)$, which is covariantly functorial (and contravariantly functorial in a special case), and homomorphisms $h_1, h_2, h_3$ such that the following diagrams commute:

$$\begin{array}{ccc} K_0(\mathcal{V}/X) & \xrightarrow{h_1} & Ab(\mathcal{V}/X) \\ & \searrow_{\chi_{Gro}} \quad \swarrow_{h_2} & \\ & K_0(\mathcal{V}) & \end{array}$$

$$\begin{array}{ccccc} K_0(\mathcal{V}/X) & \xrightarrow{h_1} & Ab(\mathcal{V}/X) & \xrightarrow{h_2} & K_0(\mathcal{V}) \\ \downarrow e & & \downarrow h_3 & & \downarrow e \\ F(X) & \xrightarrow{c_*} & H_*(X) & \xrightarrow{\int_X} & \mathbb{Z}. \end{array}$$

More on this topic will be treated in a different paper.

DEPARTMENT OF MATHEMATICS AND COMPUTER SCIENCE, FACULTY OF SCIENCE, UNIVERSITY OF KAGOSHIMA, 21-35 KORIMOTO 1-CHOME, KAGOSHIMA 890-0065, JAPAN
*E-mail address*: yokura@sci.kagoshima-u.ac.jp